\documentclass[letterpaper, 11pt]{amsart}

\usepackage{amsmath,amsthm,amsfonts,amssymb,amscd}
\usepackage{bbm}
\usepackage{bm}
\usepackage{tikz}
\usepackage{tikz-cd}
\usepackage{appendix}
\usepackage{BOONDOX-calo}
\usepackage{standalone}
\usepackage[hidelinks]{hyperref}
\usetikzlibrary{arrows,chains,matrix,positioning,scopes}
\usepackage[letterpaper, left=3cm,right=3cm, top=3cm, bottom=3cm]{geometry}
\usepackage{adjustbox}
\usepackage{enumitem}

\newcommand{\CC}{{\mathbb C}}

\newcommand{\ZZ}{{\mathbb Z}}
\newcommand{\QQ}{{\mathbb Q}}
\newcommand{\F}{{\mathbb F}}
\newcommand{\K}{{\mathbb K}}

\newcommand{\inEnd}{{\underline{End}}}
\newcommand{\fu}{{\mathfrak{u}}}

\newcommand{\fg}{{\mathfrak{g}}}
\newcommand{\fh}{{\mathfrak{h}}}

\newcommand{\inHom}{{\underline{Hom}}}

\newcommand\dual{\raise0.9ex\hbox{$\scriptscriptstyle\vee$}}

\newcommand{\bT}{{\mathbf{T}}}

\newcommand{\bS}{{\mathbf{S}}}

\newcommand{\db}{{\bullet,\bullet}}

\theoremstyle{plain}
\newtheorem{thm}{Theorem} 
\newtheorem{prop}[thm]{Proposition}
\newtheorem{lemma}[thm]{Lemma}  
\newtheorem{cor}[thm]{Corollary}
\numberwithin{thm}{subsection}
\numberwithin{equation}{section}

\newtheorem{manualtheoreminner}{Theorem}
\newenvironment{thm'}[1]{%
  \renewcommand\themanualtheoreminner{#1}%
  \manualtheoreminner
}{\endmanualtheoreminner}

\theoremstyle{definition}
\newtheorem{defn}[thm]{Definition}

\theoremstyle{remark}
\newtheorem{rem}[thm]{Remark}

\tikzset{>=stealth}

\makeatletter
\def\@seccntformat#1{%
  \protect\textup{\protect\@secnumfont
    \ifnum\pdfstrcmp{subsection}{#1}=0 \bfseries\fi
    \csname the#1\endcsname
    \protect\@secnumpunct
  }%
}  
\makeatother

\makeatletter
\@namedef{subjclassname@2020}{%
  $2020$ Mathematics Subject Classification}
\makeatother

\begin{document}

\title[Depth one part of Tannakian groups of filtrations]{Depth one part of Tannakian groups of filtrations}
\author{Payman Eskandari}
\address{Department of Mathematics and Statistics, University of Winnipeg, Winnipeg MB, Canada }
\email{p.eskandari@uwinnipeg.ca}
\begin{abstract}
Let $(F_r M)_{r\in\ZZ}$ be a finite filtration on an object $M$ of a neutral Tannakian category $\bT$ in characteristic zero. Let $\fu(M)=\fu^F(M)$ be the Lie algebra of the subgroup of the Tannakian fundamental group of $M$ that acts trivially on the associated graded $Gr^FM$. The filtration $F_\bullet M$ induces a filtration on the internal Hom $\inHom(M,M)$, which in turn induces a filtration $F_\bullet\fu(M)$ on $\fu(M)$. This filtration on $\fu(M)$ is concentrated in negative degrees. In this paper, we give a description of the graded piece $Gr^F_{-1}\fu( M)$ in terms of the extensions $F_{r+1}M/F_{r-1}M\in Ext^1(Gr^F_{r+1}M, Gr^F_{r}M)$. In particular, these extensions determine $Gr^F_{-1}\fu(M)$. Note that here we neither assume the filtration is functorial, nor we assume that $Gr^FM$ is semisimple. The problem of studying $\fu(M)$ in this generality is motivated by the desire to understand Tannakian groups associated to a mixed motive and its realizations, including realizations for which semisimplicity of realizations of pure motives is not known and realizations that lack an interesting functorial weight filtration.

We also give two related applications. The first is an equivalent condition in the generality described above for when $\fu(M)$ coincides with its trivial upper bound $F_{-1}\inHom(M,M)$. This result generalizes earlier criteria for maximality of $\fu(M)$ obtained by various authors in special contexts or under limiting conditions. In the second application, we apply the results to the setting of a neutral Tannakian category with a functorial weight filtration $W_\bullet$. Combining the constructions of \cite{Es1} with our general maximality criterion we prove a result about the structure of the set of isomorphism classes of objects $M$ for which $Gr^WM$ is isomorphic to a given graded object $A$ and $\fu^W(M)=W_{-1}\inHom(M,M)$.
\end{abstract}
\maketitle

\section{Introduction}\label{sec: Introduction}
\subsection{Contents of the paper}
Let $\bT$ be a neutral Tannakian category (in the sense of \cite{DM82}) over a field of characteristic zero. Let $M$ be an object of $\bT$ with an increasing filtration $(F_r M)_{r\in\ZZ}$ with
\begin{equation}\label{eq23}
0=F_0M \subsetneq F_1M \subsetneq F_2M\subsetneq \cdots \subsetneq F_kM=M \hspace{.4in} (k\geq 1).
\end{equation}
By Tannakian formalism, there is a canonical object $\fu(M)$ of $\bT$ with the property that for every fiber functor $\omega$, the image of $\fu(M)$ under $\omega$ is the Lie algebra of the kernel of the natural surjection from the fundamental group (in the Tannakian sense) of $M$ to the fundamental group of 
\[Gr^FM = \bigoplus_r Gr^F_rM = \bigoplus_r \, F_rM/F_{r-1}M.\] 
This kernel is easily seen to be a unipotent group, and thus the study of the fundamental group of $M$ reduces, to a large extent, to the study of $\fu(M)$ and the fundamental group of $Gr^FM$. The goal of this paper is to give some new results about $\fu(M)$, building on and generalizing some of the results of the earlier works \cite{EM1}, \cite{EM2} of the author with K. Murty and \cite{Es1} and \cite{Es2} of the author.
\medskip\par 
A natural framework in which our results are applicable is the framework of mixed motives and their various (in particular, Hodge, $\ell$-adic, and de Rham-Betti\footnote{Following \cite[\S 7.1.6]{An04}, the de Rham-Betti realization of a motive over $\K\subset\CC$ consists of a triple $(V,W,\sigma)$ consisting of its Betti realization $V$ (a vector space over $\QQ$), its de Rham realization $W$ (a vector space over $\K$), and the comparison isomorphism $\sigma: V\otimes_\QQ \CC \rightarrow W \otimes_\K \CC$. This is the relevant realization for periods.}) realizations. The filtration $(F_\bullet M)$ in this case will typically be taken to be given, up to a relabelling of the indices, by the weight filtration on the motive or its realization. In this context, one may be tempted to make the following two assumptions from the outset:
\begin{itemize}
\item[(i)] $Gr^FM$ is a semisimple object.
\item[(ii)] $\bT$ is equipped with a functorial weight filtration $W_\bullet$ similar to the weight filtration on the category of rational mixed Hodge structures or a reasonable category of mixed motives.
\end{itemize}
Our results about $\fu(M)$ will be new even under these assumptions. However, importantly, much of the paper (almost all of it with the exception of \S \ref{sec: application to classification of motives with maximal uni rads}) is carried out in the generality of an arbitrary finitely filtered object in an arbitrary Tannakian category in characteristic zero. Working in this generality is useful even for motivic applications. Indeed, for example, the semisimplicity of the $\ell$-adic and de Rham-Betti realizations of pure motives (over suitable fields) is only known in limited cases. Moreover, no interesting functorial weight filtration is known for the category of de Rham-Betti realizations.
\medskip\par 
We now proceed to give a summary of the contents of the paper. Working with a finite filtration $(F_\bullet M)$ as in equation \eqref{eq23} on an object $M$ of a Tannakian category $\bT$ over a field of characteristic zero, the filtration $F_\bullet$ induces a filtration on the internal Hom $\inEnd(M):=\inHom(M,M)$. The object $\fu(M)$ is a canonical subobject of $\inEnd(M)$, so it is equipped with a filtration defined by
\[
F_n \fu(M) := \fu(M)\cap F_n\inEnd(M) 
\]
for every integer $n$. One trivially has $F_{-1}\fu(M)=\fu(M)$. Our primary goal in this paper is to describe the graded component $Gr^F_{-1}\fu(M)$ (this is what ``depth 1" in the title of the paper means). It is easy to see that there is a canonical embedding
\begin{equation}\label{eq32}
Gr^F_{-1}\fu(M) \hookrightarrow \bigoplus_r \inHom(Gr^F_{r+1} M,Gr^F_{r} M)
\end{equation}
(see \S \ref{sec: defn of pi} and \S \ref{sec: filtration on u(M)}). The first main result of the paper is the following characterization of $Gr^F_{-1}\fu(M)$:

\begin{thm}\label{thm: main thm}
For brevity, set
\[
V_{r,r+1} \ := \ \inHom(Gr^F_{r+1}M,Gr^F_{r}M) \hspace{.2in} (1\leq r\leq k-1)
\]
and
\[
V \ := \ \bigoplus_r V_{r,r+1} \, .
\]
Let $\mathcal{E}_r$ be the element of $Ext^1_\bT(\mathbbm{1}, V_{r,r+1})$ corresponding to the extension $F_{r+1}M/F_{r-1}M$ of $Gr^F_{r+1}M$ by $Gr^F_{r}M$. Set $\mathcal{E}:=(\mathcal{E}_r)$, considered as an element of $Ext^1_\bT(\mathbbm{1},V)$. Then $Gr^F_{-1}\fu(M)$ is the smallest subobject of $V$ with the property that the pushforward of $\mathcal{E}$ along the quotient map $V\rightarrow V/Gr^F_{-1}\fu(M)$ is an extension in the Tannakian subcategory $\langle Gr^FM\rangle^{\otimes}$ of $\bT$ generated by $Gr^FM$.
\end{thm}

The result has the following immediate corollary:
\begin{cor}\label{cor: semisimple case} Let $V$ and $\mathcal{E}$ be as in Theorem \ref{thm: main thm}. 
\begin{itemize}[wide]
\item[(a)] The subobject $Gr^F_{-1}\fu(M)$ of $V$ only depends on the extensions $F_{r+1}M/F_{r-1}M$ of $Gr^F_{r+1}M$ by $Gr^F_{r}M$ (for $r=1,\ldots, k-1$).
\item[(b)] If $Gr^FM$ is semisimple, then $Gr^F_{-1}\fu(M)$ is the smallest subobject of $V$ with the property that the pushforward of $\mathcal{E}$ along the quotient map $V\rightarrow V/Gr^F_{-1}\fu(M)$ splits.
\end{itemize}
\end{cor}

The case $k=2$ of Theorem \ref{thm: main thm} (for which $\fu(M)=Gr^F_{-1}\fu(M)$) was proved by Bertrand \cite[Theorem 1.1]{Ber01} and Hardouin \cite[Th\'{e}or\`{e}me 2.1]{Har06} when $Gr^FM$ is semisimple and $Gr^F_2M=\mathbbm{1}$, and by the author and K. Murty in \cite[Theorem 3.3.1]{EM1} in general. The case $k=3$ of Theorem \ref{thm: main thm} was proved in \cite[Theorem 3.7.1]{Es2} assuming $Gr^FM$ is semisimple. The argument for Theorem \ref{thm: main thm} refines the cohomological method of the proofs of the last two aforementioned results from \cite{EM1} and \cite{Es2}. A sketch of this argument is included in \S \ref{sec: sketch of proof of main thm} below.
\medskip\par 
The paper includes some applications of Theorem \ref{thm: main thm}. In Proposition \ref{prop: maximality criteria 1} we give a necessary and sufficient criterion for when $\fu(M)$ is maximal (in the sense that it is equal to its trivial upper bound $F_{-1}\inEnd(M)$). This criterion, obtained in the fullest possible generality of finite filtrations in Tannakian categories in characteristic 0, implies in particular that whether or not $\fu(M)$ is maximal depends only on the extensions $F_{r+1}M/F_{r-1}M$ of $Gr^F_{r+1}M$ by $Gr^F_{r}M$. The result generalizes several earlier maximality criteria by various authors in special cases (in particular, Theorem 2.1 of \cite{Ber01}, Corollary 4.5 and 4.6 of \cite{BP1}, Theorem 5.3.2 of \cite{Es1} and Corollary 3.8.1 of \cite{Es2}).
\medskip\par 
In the last section of the paper we assume that the category $\bT$ has a functorial weight filtration $W_\bullet$. Given an object $M$ of $\bT$, let $\fu(M)$ be the object associated to the Lie algebra of the kernel of the restriction map from the fundamental group of $M$ to the fundamental group of $Gr^WM$. Thus if the nonzero graded components of $Gr^WM$ are in weights $p_1<\cdots<p_k$, then $\fu(M)$ coincides with its namesake from before for the filtration $F_\bullet M$ given by $F_rM:=W_{p_r}M$. Fix a graded object $A$ with $k$ graded components. In \cite[\S 3 and \S 4]{Es1} we developed a theory of ``generalized extensions" of $A$ which was used to study the set of isomorphism classes of objects $M$ of $\bT$ with associated graded $Gr^WM$ isomorphic to $A$. The case $k=3$ of this is exactly the theory of blended extensions ({\it extensions panach\'{e}es} \cite{Gr68}). This machinery was used in \cite[\S 5]{Es1} to study the set of isomorphism classes of objects $M$ of $\bT$ with $Gr^W(M)$ isomorphic to $A$ and a maximal $\fu(M)$ (i.e. with $\fu(M)=W_{-1}\inEnd(M)$). There, we were only able to give a description of this set in the case where $A$ is semisimple and ``graded-independent"; the latter hypothesis, which can be quite limiting, means that the components in a certain decomposition of $W_{-1}\inEnd(A)$ have no nonzero isomorphic subobjects (see Definition 5.3.1 of \cite{Es1}). The reason for not being able to go beyond these hypotheses was a lack of a suitable maximality criterion for $\fu(M)$ in general. With Proposition \ref{prop: maximality criteria 1} in hand, in \S \ref{sec: application to classification of motives with maximal uni rads} of this article we handle the general case with no condition on $A$. The result in recorded as Theorem \ref{thm: isomorphism classes of objects with given Gr and maximal u}.

\subsection{Sketch of the proof of Theorem \ref{thm: main thm}}\label{sec: sketch of proof of main thm}
Let $\omega$ be a fiber functor for $\bT$ and $\mathcal{G}(M)$ (respectively, $\mathcal{G}(Gr^FM)$ and $\mathcal{U}(M)$) the fundamental group of $M$ with respect to $\omega$ (respectively, the fundamental group of $Gr^FM$ with respect to $\omega$ and the kernel of the restriction map $\mathcal{G}(M)\twoheadrightarrow \mathcal{G}(Gr^FM)$). For every object $X$ of $\langle Gr^FM \rangle^{\otimes}$ we have a commutative diagram
\begin{equation}\label{eq22}
\begin{tikzcd}
0 \ar[r] & H^1(\displaystyle{\frac{\mathcal{G}(M)}{\mathcal{U}(M)}}, \omega X) \ar[r] & H^1(\mathcal{G}(M), \omega X)  \ar[r] & H^1(\mathcal{U}(M), \omega X)^{\mathcal{G}(M)} \\
& Ext^1_{\langle Gr^FM\rangle^{\otimes}}(\mathbbm{1},X) \ar[u, equal] \ar[r, phantom, "\subset"] & Ext^1_{\langle M\rangle^{\otimes}}(\mathbbm{1},X) \ar[u, equal] \ar[r, dashed, "\Psi_X"] &  Hom_\bT(\fu(M)^{ab}, X). \ar[u, equal, "(\ast)"]
\end{tikzcd}
\end{equation}
The top row is the beginning part of the inflation-restriction exact sequence for group cohomology of algebraic groups (the action of $\mathcal{U}(M)$ on $\omega X$ is trivial because $X$ is in $\langle Gr^FM \rangle^{\otimes}$). The first two vertical identifications are given by $\omega$, combined with the fact that the Yoneda $Ext^1$ groups for the category of finite-dimensional representations of an algebraic group agree with the $Ext^1$ groups defined using resolutions (i.e., with group cohomology). The identification $(\ast)$ is obtained starting from the top by noting that since the action of $\mathcal{U}(M)$ on $\omega X$ is trivial, $H^1(\mathcal{U}(M),\omega X)$ is simply the Hom group $Hom(\mathcal{U}(M),\omega X)$ in the category of algebraic groups, then passing on to the Lie algebras (since $\mathcal{U}(M)$ is unipotent and $\omega X$ is a vector group), and finally again using $\omega$ (the superscript $ab$ stands for abelianization). The map $\Psi_X$ is defined by the commutativity of the diagram. Everything is functorial in $X$ and in fact, even though not needed for the present argument, one can show that the map $\Psi_X$ is independent of the choice of $\omega$. We compute the image of the extension $\mathcal{E}$ under $\Psi_V$, where $\mathcal{E}$ and $V$ are as in the statement of Theorem \ref{thm: main thm}. We see that this image is the map induced by the natural map
\[
\fu(M) \twoheadrightarrow Gr^F_{-1}\fu(M) \subset V
\]
(see Proposition \ref{prop: main prop}). Theorem \ref{thm: main thm} is then easily deduced using functoriality of the maps $\Psi_X$ in $X$ and the fact that the kernel of $\Psi_X$ is the subgroup of extensions in $\langle Gr^FM\rangle^{\otimes}$.

\subsection{Outline of the paper}
Section \ref{sec: remarks on Ext groups} is devoted to studying the map $\Psi_X$ of equation \eqref{eq22} (in particular, its explicit description and its independence of the fiber functor). For potential future use, the section is written in a slight more generality with no extra work. In \S \ref{sec: setup and initial considerations} we introduce the setup of our problem of interest and make some initial observations. In \S \ref{sec: determination of pi(u(M))} we prove Theorem \ref{thm: main thm} (\S \ref{sec: Psi_X} - \S \ref{sec: proof of main thm}). We also give some variants of this result (\S \ref{sec: variants I} and \S \ref{sec: variants II}) and give some criteria for maximality of $\fu(M)$ and $Gr^F_{-1}\fu(M)$ (\S \ref{sec: maximality criteria} and \S \ref{sec: variants II}). Finally, in \S \ref{sec: application to classification of motives with maximal uni rads} we prove our result about the structure of the set of isomorphism classes of objects $M$ of a {\it filtered} Tannakian category with a given associated graded and maximal $\fu(M)$.

\section{$Ext^1$ groups in Tannakian categories}\label{sec: remarks on Ext groups}
The goal of this section is to make some remarks about $Ext^1$ groups in neutral Tannakian categories. In the process, we also recall some needed background and introduce some notation.
\medskip\par 
Throughout the paper, by a Tannakian category we always mean a neutral Tannakian category. The fundamentals of Tannakian formalism to the extent of \cite{DM82} will be taken for granted. For any Lie algebra object $\fg$ in any Tannakian category, we denote the abelianization $\fg/[\fg,\fg]$ by $\fg^{ab}$. Internal Homs are denoted by $\inHom$. The symbol $\mathbbm{1}$ refers to the unity object with respect to tensor product. By a Tannakian subcategory of a Tannakian category we always mean a Tannakian subcategory that is closed under taking subobjects (and hence subquotients). For any object $X$ of a Tannakian category, the Tannakian subcategory generated by $X$ is denoted by $\langle X\rangle^{\otimes}$; by definition, this is the smallest full Tannakian subcategory containing $X$.
\medskip\par 
Given a field $\K$, the category of algebraic groups over $\K$ is denoted by $\mathbf{AlgGr}(\K)$. For an object $\mathcal{G}$ of $\mathbf{AlgGr}(\K)$, the category of finite-dimensional representations of $\mathcal{G}$ over $\K$ is denoted by $\mathbf{Rep}(\mathcal{G})$.
\medskip\par 
The notation $Ext^i$ always refers to the Ext groups (or vector spaces, if applicable) in the sense of Yoneda. We include the intended category for Ext and Hom groups as subscripts (e.g., $Ext^1_\bT$ for a category $\bT$). Given an algebraic group $\mathcal{G}$ over a field $\K$, for brevity, the Hom and Ext groups in $\mathbf{Rep}(\mathcal{G})$ are simply denoted by $Hom_\mathcal{G}$ and $Ext^i_\mathcal{G}$. As usual, $Hom_\K$ means Hom in the category of vector spaces over $\K$.
\medskip\par 
Throughout the paper from this point on, we fix a field $\F$ of characteristic zero. We will simply write $\mathbf{AlgGr}$ instead of $\mathbf{AlgGr}(\F)$.

\subsection{}\label{sec: map from Ext(1,-) to Hom(u^ab,-) for alg groups} Let $\mathcal{G}$ be an algebraic group over $\F$ and $\mathcal{H}$ a normal subgroup of $\mathcal{G}$. We consider $\mathbf{Rep}(\mathcal{G}/\mathcal{H})$ as a full subcategory of $\mathbf{Rep}(\mathcal{G})$. Let $\fg$ and $\fh$ be respectively the Lie algebras of $\mathcal{G}$ and $\mathcal{H}$, both considered as representations of $\mathcal{G}$ via the adjoint action. The abelianization $\fh^{ab}$ belongs to the subcategory $\mathbf{Rep}(\mathcal{G}/\mathcal{H})$, as $\mathcal{H}$ acts on it trivially.
\medskip\par 
Let $X$ be an object of $\mathbf{Rep}(\mathcal{G}/\mathcal{H})$. Then there exists a functorial (in $X$) map
\begin{equation}\label{eq11}
Ext^1_{\mathcal{G}}(\mathbbm{1},X) \rightarrow Hom_{\mathcal{G}}(\fh^{ab}, X),
\end{equation}
defined as follows:
\medskip\par 
\noindent \underline{Step one}: The canonical isomorphism between (Yoneda) $Ext^1_{\mathcal{G}}(\mathbbm{1},X)$ and the group cohomology $H^1(\mathcal{G}, X)$ (for algebraic groups) composed with the restriction map for the subgroup $\mathcal{H}\mathrel{\unlhd} \mathcal{G}$ from group cohomology gives a map
\[
Ext^1_{\mathcal{G}}(\mathbbm{1},X) \rightarrow H^1(\mathcal{H},X)^{\mathcal{G}} = Hom_{\mathbf{AlgGr}}(\mathcal{H}, X)^{\mathcal{G}},
\]
where $X$ is considered as an additive algebraic group, and the action of $\mathcal{G}$ on $Hom_{\mathbf{AlgGr}}(\mathcal{H}, X)$ is via its natural action on $X$ and its conjugation action on $\mathcal{H}$. In the paper we will need to work with this map explicitly, so here we shall give its description more concretely (we may take this concrete description as the definition). Given $\mathcal{Z}\in Ext^1_{\mathcal{G}}(\mathbbm{1},X)$ represented by an extension
\[
\begin{tikzcd}
0 \arrow[r] & X  \arrow[r] & Z  \arrow[r] & \mathbbm{1}  \arrow[r] & 0 ,
\end{tikzcd}
\]
choose a linear section $s$ of $Z\rightarrow \mathbbm{1}=\F$ to identify $Z\cong X\oplus \F$ as vector spaces. Express the action of $\mathcal{H}$ on $Z$ in terms of this decomposition. Since $\mathcal{H}$ acts trivially on $X$, every element $h$ of $\mathcal{H}$ will act on $Z$ through an element of $GL(X\oplus \F)$ of the form
\[
\begin{pmatrix}
1 & h_{12}\\
0& 1\end{pmatrix}.
\]
Sending $h\mapsto h_{12}$ we obtain a morphism of algebraic groups
\[
\raisebox{3pt}{$\chi$}_\mathcal{Z}: \mathcal{H} \rightarrow Hom_\F(\F, X) \cong X.
\]
This map is easily seen to be independent of both the choice of the section $s$ and the choice of the representative of $\mathcal{Z}$. Thus so far, we have defined a map from $Ext^1_{\mathcal{G}}(\mathbbm{1},X)$ to $Hom_{\mathbf{AlgGr}}(\mathcal{H}, X)$.

Identifying $Z$ as $X\oplus \F$ via a choice of a section and writing elements of $GL(Z)$ as $2\times 2$ matrices again, given $g\in \mathcal{G}$ and $h\in \mathcal{H}$, the actions of the element $ghg^{-1}$ of $\mathcal{G}$ on $Z$ is through the element
\[
\begin{pmatrix}
g_{11} & g_{12}\\
0& 1\end{pmatrix} \begin{pmatrix}
1 & h_{12}\\
0& 1\end{pmatrix} \begin{pmatrix}
g_{11}^{-1} & -g_{11}^{-1}g_{12} \\
0& 1\end{pmatrix} = \begin{pmatrix}
1 & g_{11}h_{12}\\
0& 1\end{pmatrix},
\]
where $g_{11}=g_X$ is the action of $g$ on $X$. It follows that $\raisebox{3pt}{$\chi$}_{\mathcal{Z}}(ghg^{-1})= g_X\chi_{\mathcal{Z}}(h)$, so that $\raisebox{3pt}{$\chi$}_{\mathcal{Z}}$ is $\mathcal{G}$-equivariant. We have thus indeed defined a map 
\begin{equation}\label{eq12}
Ext^1_{\mathcal{G}}(\mathbbm{1},X) \rightarrow Hom_{\mathbf{AlgGr}}(\mathcal{H}, X)^{\mathcal{G}}.
\end{equation}
Its functoriality in $X$ can be checked by a direct computation from the construction.
\medskip\par 
\noindent \underline{Step two}: The map \eqref{eq11} is obtained by composing \eqref{eq12} with
\[
Hom_{\mathbf{AlgGr}}(\mathcal{H}, X)^{\mathcal{G}} \xrightarrow{Lie} Hom_{\mathbf{Lie}}(\fh, X)^{\mathcal{G}} \stackrel{(\ast)}{\cong} Hom_{\F}(\fh^{ab}, X)^{\mathcal{G}} = Hom_{\mathcal{G}}(\fh^{ab}, X),
\]
where the first map passes on to the Lie algebras ($\mathbf{Lie}$ = the category of Lie algebras over $\F$) and $(\ast)$ is because the Lie algebra $X$ is abelian.\footnote{Throughout the paper, the symbol $\cong$ will be used to refer to a canonical or distinguished isomorphism. The symbol $\simeq$ will be used to more generally mean ``isomorphism" (distinguished or otherwise) or to mean ``isomorphic" (possibly without a particular choice of an isomorphism).}

\subsection{}\label{sec: review of the very basics of tan cats} 
From this point on in the paper we fix a Tannakian category $\bT$ over $\F$. By a fiber functor we shall mean a fiber functor over $\F$. Let $\omega$ be a fiber functor for $\bT$. For any object $X$ of $\bT$, we denote the Tannakian group of $X$ with respect to $\omega$ (i.e., the group $\underline{Aut}^{\otimes}(\omega |_{\langle X\rangle^{\otimes}})$ of tensor automorphisms of the restriction of $\omega$ to $\langle X\rangle^{\otimes}$, see \cite{DM82}) by $\mathcal{G}(X,\omega)$. By the main theorem of Tannakian formalism (\cite[Theorem 2.11]{DM82}), the functor $\omega$ gives an equivalence of tensor categories
\[
\langle X\rangle^{\otimes} \rightarrow \mathbf{Rep}(\mathcal{G}(X,\omega)).
\]
We identify $\mathcal{G}(X,\omega)$ as a subgroup of $GL(\omega X)$ via its action on $X$. Via the equivalence above, the internal Hom $\inEnd(X)=\inHom(X,X)$ then becomes $End_\F(\omega X)$ equipped with the restriction of the adjoint action of $GL(\omega X)$ to the subgroup $\mathcal{G}(X,\omega)$.
\medskip\par 
From here until the end of \S \ref{sec: remarks on Ext groups}, we fix an object $M$ of $\bT$ and an object $N$ of $\langle M\rangle^{\otimes}$. There is a surjective restriction map
\[
\mathcal{G}(M,\omega) \rightarrow \mathcal{G}(N,\omega),
\]
the kernel of which we shall denote by $\mathcal{H}(M,N,\omega)$. An object $X$ of $\langle M\rangle^{\otimes}$ belongs to the subcategory $\langle N\rangle^{\otimes}$ if and only if the action of $\mathcal{H}(M,N,\omega)$ on $\omega X$ is trivial.
\medskip\par 
We shall denote the Lie algebra of $\mathcal{H}(M,N,\omega)$ by $\fh(M,N,\omega)$. It is a Lie subalgebra of $End_\F(\omega X)$. The adjoint action of $\mathcal{G}(M,\omega)$ restricts to an action on $\fh(M,N,\omega)$, which via the main theorem of Tannakian formalism gives rise to a Lie subobject
\[
\fh(M,N) \subset \inEnd(M)
\]
whose image under $\omega$ is the Lie subalgebra $\fh(M,N,\omega)$ of $End_\F(\omega X)$. Moreover, the subobject $\fh(M,N)$ of $\inEnd(M)$ is independent of the choice of the fiber functor $\omega$. See for instance, \cite[\S 2]{EM2}. 

\subsection{}\label{sec: the map Psi_X} Let $X$ be an object of the subcategory $\langle N\rangle^{\otimes}$ of $\langle M\rangle^{\otimes}$. Choose a fiber functor $\omega$ for $\bT$. In view of the equivalence of categories $\langle M\rangle^{\otimes} \xrightarrow{ \, \omega \, } \mathbf{Rep}(\mathcal{G}(M,\omega))$, the map \eqref{eq11} for $\mathcal{G}=\mathcal{G}(M,\omega)$ and $\mathcal{H}=\mathcal{H}(M,N,\omega)$ gives us a map
\[
\Psi^{M,N}_X: Ext^1_{\langle M\rangle^{\otimes}}(\mathbbm{1}, X) \rightarrow Hom_\bT(\fh(M,N)^{ab},X)
\]
that is functorial in the object $X$ of $\langle N\rangle^{\otimes}$. Note that the object $\fh(M,N)^{ab}$ belongs to $\langle N\rangle^{\otimes}$ because $\mathcal{H}(M,N,\omega)$ acts on $\omega \fh(M,N)^{ab}$ trivially.

\begin{prop}
The map $\Psi^{M,N}_X$ is independent of the choice of fiber functor $\omega$.
\end{prop}
\begin{proof}
This is a consequence of the fact, due to Deligne \cite[\S 1.12 \& \S 1.13]{De90}, that any two fiber functors are isomorphic over an algebraic closure $\overline{\F}$ of $\F$. The proof was given in details in \cite{Es2} (see Lemma 3.2.1 therein) when $\fh(M,N)$ is the Lie algebra of the unipotent radical of the fundamental group of $M$. The argument goes through essentially identically to establish the more general statement given here. Given two fiber functors $\omega$ and $\omega'$, let $\alpha$ be an isomorphism of $\otimes$-functors $\overline{\omega}\rightarrow \overline{\omega'}$, where throughout, we use the bar symbol to refer to the extension of scalars from $\F$ to $\overline{\F}$. For brevity, let us tentatively write $\fh$ for $\fh(M,N)$, and $\mathcal{G}$, $\mathcal{H}$ (resp. $\mathcal{G}'$, $\mathcal{H}'$) for $\mathcal{G}(M,\omega)$, $\mathcal{H}(M,N,\omega)$ (resp. the counterparts for $\omega'$). Then $\alpha$ gives an isomorphism $\overline{\mathcal{G}}\rightarrow \overline{\mathcal{G}'}$ restricting to an isomorphism $\overline{\mathcal{H}}\rightarrow \overline{\mathcal{H}'}$. Then there is a commutative diagram
\[
\begin{tikzcd}
& Ext^1_{\mathcal{G}}(\mathbbm{1}, \omega X) \arrow[d, hookrightarrow] \arrow[r] & Hom_{\mathcal{G}}(\omega \fh^{ab}, \omega X) \arrow[d, hookrightarrow] & \\
& Ext^1_{\overline{\mathcal{G}}}(\mathbbm{1}, \overline{\omega X})\arrow[dd, "\simeq", "\alpha" '] \arrow[r] & Hom_{\overline{\mathcal{G}}}(\overline{\omega \fh^{ab}}, \overline{\omega X}) \arrow[dd, "\simeq"', "\alpha" ] & \\
Ext^1_{\langle M\rangle^{\otimes}}(\mathbbm{1}, X) \arrow[ruu, "\text{$\omega, \simeq$}"] \arrow[rdd, "\text{$\omega', \simeq$}" '] & & & Hom_\bT(\fh^{ab},X). \arrow[luu, "\text{$\omega, \simeq$}" '] \arrow[ldd, "\text{$\omega', \simeq$}"] \\
& Ext^1_{\overline{\mathcal{G}'}}(\mathbbm{1}, \overline{\omega' X}) \arrow[r] & Hom_{\overline{\mathcal{G}'}}(\overline{\omega' \fh^{ab}}, \overline{\omega' X}) & \\
& Ext^1_{\mathcal{G}'}(\mathbbm{1}, \omega' X) \arrow[u, hookrightarrow] \arrow[r] & Hom_{\mathcal{G}'}(\omega' \fh^{ab}, \omega' X) \arrow[u, hookrightarrow] & 
\end{tikzcd}
\]
Here, the horizontal maps are all of the form of equation \eqref{eq11} of \S \ref{sec: map from Ext(1,-) to Hom(u^ab,-) for alg groups}. The vertical maps in the top and bottom squares are by extension of scalars; that the ones on the left (from $Ext^1_{\mathcal{G}}$ to $Ext^1_{\overline{\mathcal{G}}}$ and the counterpart for $\mathcal{G}'$) are injective will not matter for the proof, but it is by the fact that invariants behave well with respect to extension of scalars \cite[Proposition 4.31]{Mi17}. The isomorphisms marked by $\omega$ and $\omega'$ are given by Tannakian formalism. The isomorphism $\alpha$ on the right is given by the isomorphism of functors $\alpha$. The isomorphism marked by $\alpha$ on the left is obtained at the level of extensions by replacing the $\overline{\mathcal{G}}$-representation $\overline{\omega X}$ by the $\overline{\mathcal{G}'}$-representation $\overline{\omega' X}$ via $\alpha$, keeping the middle vector space unchanged, and transferring the $\overline{\mathcal{G}}$-action on it to a $\overline{\mathcal{G}'}$-action via $\alpha$. That the triangles on the left and right commute is because $\alpha$ is an isomorphism of functors. That the top and bottom squares commute is easily seen from the definition of the map \eqref{eq11}. That the middle square is commutative is seen via again the definition of the horizontal maps, on noting that the map $\overline{\omega \fh}\rightarrow \overline{\omega' \fh}$ given by the morphism of functors $\alpha$ coincides with the map obtained by first considering the isomorphism $\overline{\mathcal{H}}\rightarrow \overline{\mathcal{H}'}$ given by $\alpha$, and then passing to the Lie algebras.

The independence of $\Psi^{M,N}_X$ from the choice of $\omega$ now follows from the commutativity of the diagram and injectivity of the top and bottom vertical arrows on the right. 
\end{proof}

\subsection{} \label{sec: summary of construction of map Ext^1(1,-) -> Hom(h^ab,-)}
For convenience and referencing purposes, let us summarize the discussion so far. Fix an object $M$ of $\bT$ ( = a Tannakian category over $\F$) and an object $N$ of $\langle M\rangle^{\otimes}$. Then for every object $X$ of $\langle N\rangle^{\otimes}$ we have constructed a map
\begin{equation}\label{eq13}
\Psi^{M,N}_X: Ext^1_{\langle M\rangle^{\otimes}}(\mathbbm{1}, X) \rightarrow Hom_\bT(\fh(M,N)^{ab},X),
\end{equation}
which is functorial with respect to morphisms $X\rightarrow X'$ between two objects $X,X'$ of $\langle N\rangle^{\otimes}$. Given an extension class $\mathcal{Z}$ of $\mathbbm{1}$ by $X$ in $\langle M\rangle^{\otimes}$, the image of $\mathcal{Z}$ under $\Psi^{M,N}_X$ is computed as follows: Choose a representative extension for $\mathcal{Z}$ with its middle denoted by $Z$. Choose a fiber functor $\omega$ and a linear section of the surjection $\omega Z\rightarrow \omega\mathbbm{1}=\F$ to obtain a decomposition $\omega Z\cong \omega X\oplus \F$ as vector spaces. The action of $\mathcal{H}(M,N,\omega)$ on $\omega Z$ expressed in terms of this decomposition gives rise to a $\mathcal{G}(M,\omega)$-equivariant morphism of algebraic groups $\mathcal{H}(M,N,\omega)\rightarrow \omega X$. Passing to the Lie algebras and the abelianization, we obtain a morphism of $\mathcal{G}(M,\omega)$-representations $\omega \fh(M,N)^{ab} = \fh(M,N,\omega)^{ab}\rightarrow \omega X$. This map is $\omega$ of $\Psi^{M,N}_X(\mathcal{Z})$. This recipe for defining $\Psi^{M,N}_X(\mathcal{Z})$ is indeed well-defined and in particular, is independent of the choice of a fiber functor $\omega$.

\subsection{}\label{sec: kernel of Psi_X}
Under mild conditions (which will be satisfied for the purposes of this paper), the kernel of the map $\Psi^{M,N}_X$ can be intrinsically described, as follows. 

\begin{lemma}\label{lem: kernel of Psi_X}
As in \S \ref{sec: the map Psi_X} and \S \ref{sec: summary of construction of map Ext^1(1,-) -> Hom(h^ab,-)}, let $N$ be an object of $\langle M\rangle^{\otimes}$ and $X$ an object of $\langle N\rangle^{\otimes}$. Suppose that the subgroup $\mathcal{H}(M,N,\omega)$ of $\mathcal{G}(M,\omega)$ is connected for a choice of fiber functor $\omega$. Then the kernel of $\Psi^{M,N}_X$ is equal to the subgroup $Ext^1_{\langle N\rangle^{\otimes}}(\mathbbm{1}, X)$ of $Ext^1_{\langle M\rangle^{\otimes}}(\mathbbm{1}, X)$.
\end{lemma}

\begin{proof}
We calculate $\Psi^{M,N}_X$ using the fiber functor $\omega$. By construction (see \S \ref{sec: summary of construction of map Ext^1(1,-) -> Hom(h^ab,-)} and \S \ref{sec: map from Ext(1,-) to Hom(u^ab,-) for alg groups}), $\Psi^{M,N}_X$ is the following composition:
\[
\begin{tikzcd}[column sep = small]
Ext^1_{\langle M\rangle^{\otimes}}(\mathbbm{1}, X) \arrow[r, "\omega", "\simeq" '] & Ext^1_{\mathcal{G}(M,\omega)}(\mathbbm{1}, \omega X) \xrightarrow[(\ast)]{\text{Eq.\eqref{eq12}}} Hom_{\mathbf{AlgGr}}(\mathcal{H}(M,N,\omega), \omega X)^{\mathcal{G}(M,\omega)} \arrow[d, shift left=23, "Lie"]\\
Hom_\bT(\fh(M,N)^{ab}, X) & \arrow[l, "\omega" ', "\simeq"] Hom_{\mathcal{G}(M,\omega)}(\fh(M,N,\omega)^{ab}, \omega X) \cong Hom_{\mathbf{Lie}}(\fh(M,N,\omega), \omega X)^{\mathcal{G}(M,\omega)}
\end{tikzcd}
\]
By the hypothesis of connectedness of $\mathcal{H}(M,N,\omega)$, the map denoted by $Lie$ (which passes to the Lie algebras) is injective. Thus an extension $\mathcal{Z}\in Ext^1_{\langle M\rangle^{\otimes}}(\mathbbm{1}, X)$ is in the kernel of $\Psi^{M,N}_X$ if and only if $\omega \mathcal{Z}$ is in the kernel of the map $(\ast)$ of the diagram. The result now follows from the first three terms of the inflation-restriction exact sequence for cohomology of algebraic groups together with the interpretation of $H^1(\mathcal{G}(M,\omega)/\mathcal{H}(M,N,\omega), \omega X)$ as $Ext^1_{\langle N\rangle^{\otimes}}(\mathbbm{1}, X)$. In the interest of completeness and since we already have an explicit definition of $(\ast)$ at hand, we give a more explicit version of the argument: Denoting the middle object of a representative of $\mathcal{Z}$ by $Z$, upon recalling the construction of the map \eqref{eq12} from \S \ref{sec: map from Ext(1,-) to Hom(u^ab,-) for alg groups}, we see that $\Psi^{M,N}_X(\mathcal{Z})=0$ if and only if $\omega \mathcal{Z}$ admits an $\mathcal{H}(M,N,\omega)$-equivariant section, which in turn (since $X\in \langle N\rangle^{\otimes}$) is equivalent to the action of $\mathcal{H}(M,N,\omega)$ on $\omega Z$ being trivial. The latter is equivalent to $Z$ being in $\langle N\rangle^{\otimes}$, or in other words, $\mathcal{Z}$ being in $Ext^1_{\langle N\rangle^{\otimes}}(\mathbbm{1}, X)$.
\end{proof}

\subsection{}\label{sec: when h(M) is the Lie algebra of the unipotent radical} We now consider an important special case. Suppose that $\langle N\rangle^{\otimes}$ is the full subcategory of all the semisimple objects of $\langle M\rangle^{\otimes}$, or that equivalently, $\mathcal{H}(M,N,\omega)$ for one or all choices of $\omega$ is the unipotent radical of $\mathcal{G}(M,\omega)$. Then for every semisimple object $X$ of $\langle M\rangle^{\otimes}$, the map $\Psi^{M,N}_X$ of equation \eqref{eq13} is an isomorphism. Indeed, the injectivity of $\Psi^{M,N}_X$ follows from Lemma \ref{lem: kernel of Psi_X}, as $\mathcal{H}(M,N,\omega)$ is connected and $\langle N\rangle^{\otimes}$ is semisimple. As for surjectivity, since $\mathcal{H}(M,N,\omega)$ is unipotent, the map $Lie$ of the diagram of the proof of Lemma \ref{lem: kernel of Psi_X} is an isomorphism, so the surjectivity amounts to the surjectivity of the map marked as $(\ast)$ in the same diagram. The surjectivity of this map can be seen via the inflation-restriction sequence of group cohomology (since $\mathcal{G}(M,\omega)/\mathcal{H}(M,N,\omega)\cong \mathcal{G}(N,\omega)$ is reductive and hence $H^2(\mathcal{G}(N,\omega),-)$ vanishes). One can also easily check surjectivity of $(\ast)$ in this case just from our explicit definition from \S \ref{sec: map from Ext(1,-) to Hom(u^ab,-) for alg groups}. See \cite[\S 3.1]{Es2} for details.

\section{The setup and initial considerations}\label{sec: setup and initial considerations}
Recall that $\bT$ is a (neutral) Tannakian category over a field $\F$ of characteristic 0.
\subsection{}\label{sec: introducing the basic setting}
From here until the end of \S \ref{sec: determination of pi(u(M))}, we fix the following data:
\begin{itemize}
\item[-] an object $M$ of $\bT$, equipped with a finite increasing filtration $(F_r M)_{r\in\ZZ}$ with
\[
0=F_0M \subsetneq F_1M \subsetneq F_2M\subsetneq \cdots \subsetneq F_kM=M \hspace{.4in}(k\geq 1).
\]
\end{itemize}
As usual, we denote
\[
Gr^FM = \bigoplus\limits_r Gr^F_r M = \bigoplus\limits_r \frac{F_rM}{F_{r-1}M}.
\]
We will adopt the notation of \S \ref{sec: review of the very basics of tan cats} for Tannakian fundamental groups, their subgroups, Lie algebras, and etc. Throughout, we shall set
\[
\fu(M) := \fh(M,Gr^FM)
\]
and for every fiber functor $\omega$,
\[
\mathcal{U}(M,\omega) := \mathcal{H}(M,Gr^FM,\omega)
\]
(see \S \ref{sec: review of the very basics of tan cats}). The fundamental group $\mathcal{G}(M,\omega)$ of $M$ respects the filtration $F_\bullet \omega M:=\omega F_\bullet M$, and $\mathcal{U}(M,\omega)$ is the subgroup of $\mathcal{G}(M,\omega)$ consisting of the elements that act trivially on $Gr^F\omega M$. In particular, $\mathcal{U}(M,\omega)$ is unipotent. If $Gr^FM$ is semisimple, then $\mathcal{U}(M,\omega)$ will be the unipotent radical of $\mathcal{G}(G,\omega)$. In what follows, unless otherwise indicated, we will {\it not} assume that $Gr^FM$ is semisimple.

\subsection{} The filtration $F_\bullet$ on $M$ induces a finite increasing filtration on the internal Hom $\inEnd(M)$. For any fiber functor $\omega$ and any integer $n$, the image of $F_n\inEnd(M)$ under $\omega$ is the subspace of $\omega \inEnd(M) = End_\F(\omega M)$ consisting of all the linear maps $f: \omega M\rightarrow \omega M$ such that $f(F_r \omega M) \subset F_{r+n}\omega M$ for all $r$. That the subobject $F_n\inEnd(M)$ defined this way is independent of the choice of $\omega$ can be seen using the fact, due to Deligne \cite{De90}, that every two fiber functors are isomorphic after base change to an algebraic closure of $\F$. From the description of $F_\bullet\inEnd(M)$ it is easily seen that
\[
[F_n\inEnd(M), F_m\inEnd(M)] = F_{n+m}\inEnd(M).
\]
In particular, $F_{-1}\inEnd(M)$ is a Lie subobject of $\inEnd(M)$ with derived algebra $F_{-2}\inEnd(M)$.

\subsection{}\label{sec: splitting} Let $\omega$ be a fiber functor. By a {\it splitting of $\omega M$} we shall mean the data of a (linear) section $s_r$ of the quotient map $F_r\omega M \rightarrow Gr^F_r\omega M$ for each integer $r$. Given a choice of a splitting  $(s_r)$ of $\omega M$, for each $r$ we use the section $s_r$ to obtain an isomorphism $F_r \omega M \cong F_{r-1}\omega M \oplus Gr^F_r\omega M$ as vector spaces. Putting these isomorphisms together we obtain isomorphisms
\[
\omega M = F_k \omega M \cong F_{k-1} \omega M\oplus Gr^F_k\omega M \cong F_{k-2} \omega M\oplus Gr^F_{k-1}\omega M\oplus Gr^F_k\omega M \cong \cdots
\]
and thus an isomorphism
\begin{equation}\label{eq3}
\omega M \cong \bigoplus\limits_{r=1}^k Gr^F_r\omega M
\end{equation}
of vector spaces. This isomorphism induces isomorphisms
\begin{equation}\label{eq9}
\frac{F_n\omega M}{F_m\omega M} \cong \bigoplus\limits_{r=m+1}^n Gr^F_r\omega M.
\end{equation}
for all $m<n$, all compatible with the natural inclusion and projection maps.

Once a splitting of $\omega M$ is chosen, we make identifications given by \eqref{eq3} and \eqref{eq9}. We will then write elements of $End_\F(\omega M)$ and $GL(\omega M)$ as $k\times k$ matrices. Similarly, elements of $End_\F(\frac{F_n\omega M}{F_m\omega M})$ and more generally, $Hom_\F(\frac{F_n\omega M}{F_m\omega M}, \frac{F_{n'}\omega M}{F_{m'}\omega M})$ will become matrices. The entry $f_{ij}$ of an element $f=(f_{ij})\in End_\F (\omega M)$ is the component of $f$ in $Hom_\F(Gr^F_j\omega M, Gr^F_i\omega M)$ via the decomposition
\begin{equation}\label{eq26}
End_\F \omega M \cong End_\F (Gr^F\omega M) = \bigoplus\limits_{1\leq i,j\leq k} Hom_\F(Gr^F_j\omega M, Gr^F_i\omega M).
\end{equation}
The group $\mathcal{G}(M,\omega)$ is contained in the subgroup of $GL(\omega M)$ consisting of upper triangular invertible matrices. If $\sigma=(\sigma_{ij})$ is an element of $\mathcal{G}(M,\omega)$, the automorphism $\sigma_{\frac{F_nM}{F_mM}}$ (i.e. the action of $\sigma$ on $\frac{F_nM}{F_mM}$) of $F_n\omega M/F_m\omega M$ is simply given by truncating the matrix $(\sigma_{ij})$ to the part with $m< i,j\leq n$. 

The unipotent group $\mathcal{U}(M,\omega)$ is the intersection of $\mathcal{G}(M,\omega)$ and the subgroup of $GL(\omega M)$ consisting of upper triangular elements with the identity maps on the diagonal. The Lie algebra $\omega F_{-1}\inEnd(M)$ consists of strictly upper triangular elements of $End_\F(\omega M)$. The subalgebra $\omega F_{-2}\inEnd(M)$ consists of those matrices $(f_{ij})$ for which $f_{ij}=0$ if $j-i\leq 1$, that is, it consists of strictly upper triangular matrices with zeros on the superdiagonal.

\subsection{}\label{sec: defn of pi} The goal of this paragraph is to define a canonical surjective morphism
\[
\pi: F_{-1}\inEnd(M) \rightarrow \bigoplus\limits_{r} \inHom(Gr^F_{r+1}M, Gr^F_{r}M)
\]
with
\[
\ker(\pi) = F_{-2}\inEnd(M).
\]
This map will play an important role in the paper. Let $\omega$ be a fiber functor for $\bT$. We will first explicitly define the linear map $\omega \pi$, and then verify that $\omega \pi$ indeed comes from a morphism $\pi$ that is independent of the choice of $\omega$.
\medskip\par 
Choose a splitting of $\omega M$ in the sense of \S \ref{sec: splitting}. Use the isomorphism \eqref{eq26} given by this splitting to write the elements of $End_\F(\omega M)$ as $k\times k$ matrices. Let $\omega \pi$ be the linear map 
\begin{align*}
\omega F_{-1}\inEnd(M) & = \bigoplus\limits_{\stackrel{i,j}{j-i\geq 1}} Hom_\F(Gr^F_{j}\omega M, Gr^F_{i}\omega M)  \xrightarrow{ \ \omega\pi \ } \\ & \bigoplus\limits_{r} \ Hom_\F(Gr^F_{r+1}\omega M, Gr^F_{r}\omega M) = \bigoplus\limits_{r} \ \omega \inHom(Gr^F_{r+1}M, Gr^F_{r}M)
\end{align*}
that sends
\[
f=(f_{i,j}) \mapsto (f_{r,r+1})
\]
(the indices compatible with those used in the decompositions). That is, $\omega \pi$ is the projection onto the superdiagonal entries. The notation $\omega \pi$ for this map will be justified momentarily in the next lemma. Clearly, $\omega \pi$ is surjective and 
\[
\ker(\omega\pi) = \bigoplus\limits_{\stackrel{i,j}{j-i\geq 2}} Hom_\F(Gr^F_{j}\omega M, Gr^F_{i}\omega M) = \omega F_{-2}\inEnd(M).
\]

\begin{lemma}\label{lem: defn of pi}
The map 
\[
\omega \pi: \omega F_{-1}\inEnd(M) \rightarrow \bigoplus\limits_{r} \ \omega \inHom(Gr^F_{r+1}M, Gr^F_{r}M)
\]
is independent of the choice of splitting of $\omega M$, and is the image under $\omega$ of a morphism
\[
F_{-1}\inEnd(M) \rightarrow \bigoplus\limits_{r} \inHom(Gr^F_{r+1}M, Gr^F_{r}M)
\]
denoted by $\pi$. Moreover, the morphism $\pi$ is independent of the choice of fiber functor $\omega$.
\end{lemma}

\begin{proof}
We will give an intrinsic description of $\pi$. All the claims made in the construction below can be checked easily (upon taking a fiber functor $\omega$ and a splitting of $\omega M$, if needed).

For each $r$, there is a canonical surjective map
\begin{equation}\label{eq24}
\inHom(\frac{M}{F_rM}, F_rM) \rightarrow \inHom(Gr^F_{r+1}M, Gr^F_{r}M)
\end{equation}
given by functoriality properties of $\inHom$ (say, first apply $\inHom(-,F_rM)$ to $Gr^F_{r+1}M\hookrightarrow M/F_rM$, then $\inHom(Gr^F_{r+1}M,-)$ to $F_rM\twoheadrightarrow Gr^F_{r}M$). The kernel of \eqref{eq24} is the subobject of $\inHom(\frac{M}{F_rM}, F_rM)$, denoted by $F_{-2}\inHom(\frac{M}{F_rM}, F_rM)$, that after applying a fiber functor $\omega$ consists of all linear maps $\omega M/F_r\omega M \xrightarrow {f} F_r\omega M$ satisfying $f(Gr^F_{r+1}\omega M)\subset F_{r-1}\omega M$. Taking the direct sum of the maps \eqref{eq24} over $r$ we obtain a surjective map
\begin{equation} \label{eq1}
\bigoplus\limits_r \inHom(\frac{M}{F_rM}, F_rM) \rightarrow \bigoplus\limits_r \inHom(Gr^F_{r+1}M, Gr^F_{r}M).
\end{equation}
On the other hand, for each $r$ there is an obvious canonical inclusion
\[
\inHom(\frac{M}{F_rM}, F_rM) \hookrightarrow F_{-1}\inEnd(M),
\]
again given by functoriality of $\inHom$. The summation map
\[
\sum: \  \bigoplus\limits_r \inHom(\frac{M}{F_rM}, F_rM) \rightarrow F_{-1}\inEnd(M) 
\]
is surjective. Moreover, the map \eqref{eq1} vanishes on the kernel of the summation map above. The map $\pi$ is then the map fitting in the following commutative diagram:
\[
\begin{tikzcd}[column sep = large]
\bigoplus\limits_r \inHom(\displaystyle{\frac{M}{F_rM}}, F_rM) \ar[twoheadrightarrow, rd, "\eqref{eq1}"] \ar[twoheadrightarrow, d,"\sum" '] & \\
F_{-1}\inEnd(M) \ar[r, "\pi"] & \bigoplus\limits_r \inHom(Gr^F_{r+1}M, Gr^F_{r}M).
\end{tikzcd}
\]
The image of $\pi$ under any fiber functor $\omega$ is the map denoted by $\omega \pi$ earlier.
\end{proof}
Let
\[
\pi_{r,r+1} : F_{-1}\inEnd(M) \rightarrow \inHom(Gr^F_{r+1}M, Gr^F_{r}M)
\]
be the composition of $\pi$ with the projection to $\inHom(Gr^F_{r+1}M, Gr^F_{r}M)$, so that $\pi=(\pi_{r,r+1})$. After choosing a fiber functor $\omega$ and a splitting of $\omega M$, the map $\omega \pi_{r,r+1}$ sends an element $f=(f_{ij})\in \omega F_{-1}\inEnd(M)$ to the superdiagonal entry $f_{r,r+1}$.

\subsection{}\label{sec: filtration on u(M)}
The elements of $\mathcal{U}(M,\omega)$ stabilize the filtration $F_\bullet \omega M$ and induce identity on $Gr^F\omega M$. Thus
\[
\fu(M) \subset F_{-1}\inEnd(M).
\]
Define a filtration $F_\bullet \fu(M)$ on $\fu(M)$ by setting 
\[
F_n\fu(M) := \fu(M)\cap F_n\inEnd(M).
\]
Then $Gr^F \fu(M)$ is concentrated in negative degrees. Our focus in this paper is on the graded component $Gr^F_{-1}\fu(M)$. In view of \S \ref{sec: defn of pi}, we have a commutative diagram
\begin{equation}\label{eq2}
\begin{tikzcd}
0 \arrow[r] & F_{-2}\fu(M) \arrow[r] \arrow[d, hookrightarrow] & \fu(M) \arrow[r] \arrow[d, hookrightarrow] & Gr^F_{-1} \fu(M) \arrow[r] \arrow[d, hookrightarrow] & 0\\
0 \arrow[r] & F_{-2}\inEnd(M) \arrow[r] & F_{-1}\inEnd(M) \arrow[r, "\pi"] & \bigoplus\limits_r \inHom(Gr^F_{r+1}M, Gr^F_{r}M) \arrow[r] & 0
\end{tikzcd}
\end{equation}
with exact rows. We identify $Gr^F_{-1} \fu(M)$ as a subobject of $\bigoplus\limits_r \inHom(Gr^F_{r+1}M, Gr^F_{r}M)$ via the vertical map on the right (this is the embedding in equation \eqref{eq32} of \S \ref{sec: Introduction}). Our goal is to determine 
\[
Gr^F_{-1} \fu(M) = \pi(\fu(M)) \subset \bigoplus\limits_r \inHom(Gr^F_{r+1}M, Gr^F_{r}M).
\]

\section{Determination of $Gr^F_{-1}\fu(M)$}\label{sec: determination of pi(u(M))}
We will prove Theorem \ref{thm: main thm} in this section. In addition, we will also give some variants and consequences of the result. Throughout the section, we shall work in the setting of \S \ref{sec: introducing the basic setting}, adopting the notation as introduced in the rest of \S \ref{sec: setup and initial considerations} and in \S \ref{sec: remarks on Ext groups}.  

\subsection{The maps $\Psi_X$}\label{sec: Psi_X}
Out tool to study $Gr^F_{-1}\fu(M)$ will be the maps of the form $\Psi^{M,N}_X$ constructed in \S \ref{sec: remarks on Ext groups}, taking $N=Gr^FM$. By \S \ref{sec: map from Ext(1,-) to Hom(u^ab,-) for alg groups} - \S \ref{sec: summary of construction of map Ext^1(1,-) -> Hom(h^ab,-)}, for every object $X$ of $\langle Gr^FM\rangle^{\otimes}$ we have a canonical map
\begin{equation}\label{eq4}
\Psi_X:= \Psi^{M,Gr^FM}_X: Ext^1_{\langle M\rangle^{\otimes}}(\mathbbm{1},X) \rightarrow  Hom_\bT(\fu(M)^{ab},X),
\end{equation}
which is functorial in $X$. This map can be computed by choosing any fiber functor $\omega$, following the procedure outlined in \S \ref{sec: summary of construction of map Ext^1(1,-) -> Hom(h^ab,-)} (with $\mathcal{H}(M,N,\omega)= \mathcal{U}(M,\omega)$). Since $\mathcal{U}(M,\omega)$ is unipotent, by Lemma \ref{lem: kernel of Psi_X} we have
\[
\ker(\Psi_X) = Ext^1_{\langle Gr^FM\rangle^{\otimes}}(\mathbbm{1},X).
\]
If $Gr^FM$ is semisimple, then the map $\Psi_X$ is injective for every object $X$ of $\langle Gr^F M\rangle^{\otimes}$ (in fact, it is an isomorphism because we are then in the situation of \S \ref{sec: when h(M) is the Lie algebra of the unipotent radical}).

The quotient map $\fu(M)\rightarrow \fu(M)^{ab}$ induces an injection
\[
Hom_\bT(\fu(M)^{ab}, X) \hookrightarrow Hom_\bT(\fu(M), X).
\]
For every object $X$ of $\langle Gr^FM\rangle^{\otimes}$ we use this injection to consider $\Psi_X$ as a map
\begin{equation}\label{eq5}
\Psi_X: Ext^1_{\langle M\rangle^{\otimes}}(\mathbbm{1}, X) \rightarrow Hom_\bT(\fu(M), X)
\end{equation}
(functorial in $X$ and with the same kernel as \eqref{eq4}). The image of an extension class $\mathcal{Z}$ of $\mathbbm{1}$ by $X$ in $\langle M\rangle^{\otimes}$ under \eqref{eq5} is computed using the same procedure as in \S \ref{sec: summary of construction of map Ext^1(1,-) -> Hom(h^ab,-)}, except we skip the last step of passing to the abelianization: Denoting the middle object of a representative of $\mathcal{Z}$ by $Z$, the logarithm of the morphism $\mathcal{U}(M,\omega)\rightarrow \omega X$ that describes the action of $\mathcal{U}(M,\omega)$ on $\omega Z$ is equal to the image under $\omega$ of the morphism $\fu(M)\rightarrow X$ corresponding to $\mathcal{Z}$ under \eqref{eq5}.

\subsection{The key computation}\label{sec: the main proposition} In this subsection we establish the key technical component of the proof of Theorem \ref{thm: main thm}. Following \S \ref{sec: Introduction}, we shall set
\[
V_{r,r+1}:= \inHom(Gr^F_{r+1}M, Gr^F_{r}M) \hspace{.3in}\text{and} \hspace{.2in} V:=\bigoplus\limits_r V_{r,r+1}.
\]
Recall the definitions of the extensions $\mathcal{E}_r$ and $\mathcal{E}$ from \S \ref{sec: Introduction}: The filtration $F_\bullet M$ gives rise to an extension 
\begin{equation}\label{eq7}
\begin{tikzcd}
0 \arrow[r] & Gr^F_rM \arrow[r] & \displaystyle{\frac{F_{r+1}M}{F_{r-1}M}} \arrow[r] & Gr^F_{r+1}M \arrow[r] & 0
\end{tikzcd}
\end{equation}
for each $r$. The element
\[\mathcal{E}_r \in Ext^1_{\langle M\rangle^{\otimes}}(\mathbbm{1}, V_r)\]
is the extension class corresponding to \eqref{eq7} under the canonical isomorphism
\[
Ext^1_{\langle M\rangle^{\otimes}}(\mathbbm{1}, V_{r,r+1}) \cong Ext^1_{\langle M\rangle^{\otimes}}(Gr^F_{r+1}M,Gr^F_{r}M).
\]
The extension class
\[
\mathcal{E} \in Ext^1_{\langle M\rangle^{\otimes}}(\mathbbm{1}, V)
\]
is the element corresponding to the tuple $(\mathcal{E}_r)\in \bigoplus\limits_r Ext^1_{\langle M\rangle^{\otimes}}(\mathbbm{1}, V_{r,r+1})$ under the canonical isomorphism
\[
\bigoplus_r Ext^1_{\langle M\rangle^{\otimes}}(\mathbbm{1}, V_{r,r+1}) \cong Ext^1_{\langle M\rangle^{\otimes}}(\mathbbm{1}, V).
\] 

Recall from \S \ref{sec: defn of pi} that we have a map 
\[\pi: F_{-1}\inEnd(M) \rightarrow V,\]
whose component in $V_{r,r+1}$ was denoted by $\pi_{r,r+1}$. Denote the map obtained by restricting the domain of $\pi$ (resp. $\pi_{r,r+1}$) to $\fu(M)$ also by $\pi$ (resp. $\pi_{r,r+1}$). Recall that after applying a fiber functor $\omega$ and choosing a splitting of $\omega M$, the map $\pi_{r,r+1}$ simply sends $(f_{ij})\mapsto f_{r,r+1}$.

Finally, recall that for every object $X$ of $\langle Gr^FM\rangle^{\otimes}$ we have a map $\Psi_X$ as in \eqref{eq5} of \S \ref{sec: Psi_X}.

\begin{prop}\label{prop: main prop} ~
\begin{itemize}[wide]
\item[(a)] $\Psi_{V_{r,r+1}}(\mathcal{E}_r) = \pi_{r,r+1}$ for every $r$.
\item[(b)] $\Psi_V(\mathcal{E}) = \pi$
\end{itemize}
\end{prop}

\begin{proof} \begin{itemize}[wide]
\item[(a)] Fix $r$. Let $\omega$ be a fiber functor. Recalling the explicit description of the canonical isomorphism $Ext^1(N,L)\cong Ext^1(\mathbbm{1}, \inHom(N,L))$ for any objects $N$ and $L$ of a Tannakian category (see for example, \cite[\S 3.2]{EM1}), $\mathcal{E}_r$ is given by the extension
\begin{equation}\label{eq8}
\begin{tikzcd}
0 \arrow[r] & V_{r,r+1} \arrow[r] & E_{r,r+1} \arrow[r] & \mathbbm{1} \arrow[r] & 0,
\end{tikzcd}
\end{equation}
where $E_{r,r+1}$ is the object $\inHom(Gr^F_{r+1}M, F_{r+1}M/F_{r-1}M)^\dagger$ in the notation of \cite{EM1}: This is the subobject of $\inHom(Gr^F_{r+1}M, F_{r+1}M/F_{r-1}M)$ whose image under $\omega$ consists of all linear maps $g: Gr^F_{r+1}\omega M\rightarrow F_{r+1}\omega M/F_{r-1}\omega M$ such that the composition 
\[Gr^F_{r+1}\omega M \xrightarrow{ \ g \ } \displaystyle{\frac{F_{r+1}\omega M}{F_{r-1}\omega M}}  \ \twoheadrightarrow \ Gr^F_{r+1}\omega M\]
(the latter being the canonical surjection) is a scalar multilple $\lambda(g)\cdot \text{Id}$ of the identity map. The injection $V_{r,r+1} \hookrightarrow E_{r,r+1}$ is given by functoriality of $\inHom$, and the surjection $E_{r,r+1}\twoheadrightarrow \mathbbm{1}$ after applying $\omega$ is simply the map $\lambda: g\mapsto \lambda(g)$. See \cite[\S 3.2]{EM1} for more details.

Choose a splitting of $\omega M$ in the sense of \S \ref{sec: splitting} to identify the vector space $\omega M$ with $Gr^F\omega M$. Let $s_{r+1}$ be the section of the quotient map $F_{r+1}\omega M \rightarrow Gr^F_{r+1}\omega M$ used to form the splitting. Then $s_{r+1}$ gives a section $\overline{s_{r+1}} = s_{r+1} \pmod{F_{r-1}\omega M}$ of the map $F_{r+1}\omega M/F_{r-1}\omega M \rightarrow Gr^F_{r+1}\omega M$, which in turn gives rise to a section of the map $\lambda$. We will use this section (which is given by $1\mapsto \overline{s_{r+1}}$) to compute $\Psi_{V_{r,r+1}}(\mathcal{E}_r)$ following the steps outlined at the end of \S \ref{sec: Psi_X}. 

Identify
\begin{equation}\label{eq10}
\omega E_{r,r+1} \cong Hom_\F(Gr^F_{r+1}\omega M, Gr^F_{r}\omega M) \oplus \F 
\end{equation}
via the section of $\lambda$ above. The map $\overline{s_{r+1}}$ gives an isomorphism $F_{r+1}\omega M/F_{r-1}\omega M \cong Gr^F_{r}\omega M \oplus Gr^F_{r+1}\omega M$, which agrees with the one in equation \eqref{eq9} of \S \ref{sec: splitting} for $n=r+1$ and $m=r-1$ and our choice of splitting of $\omega M$. This gives an isomorphism
\[
Hom_\F(Gr^F_{r+1}\omega M, \displaystyle{\frac{F_{r+1}\omega M}{F_{r-1}\omega M}}) \cong Hom_\F(Gr^F_{r+1}\omega M, Gr^F_{r}\omega M) \oplus Hom_\F(Gr^F_{r+1}\omega M, Gr^F_{r+1}\omega M),
\]
compatible with \eqref{eq10} via the natural inclusion maps (embedding $\F$ as the scalar endomorphisms of $Gr^F_{r+1}\omega M$). 

Use our splitting to write endomorphisms of $\omega M$ as $k\times k$ matrices. Let $f=(f_{ij})\in \omega\fu(M)$ (a strictly upper triangular matrix). Set $\sigma:=\exp(f)\in \mathcal{U}(M,\omega)$. Then $f$ and $\sigma$ have the same superdiagonal entries because $f^n\in \omega F_{-2}\inEnd(M)$ for $n\geq 2$. The action of $\sigma$ on $\omega\inHom(Gr^F_{r+1}M, F_{r+1}M/ F_{r-1}M)$ is by sending $g:Gr^F_{r+1}\omega M \rightarrow F_{r+1}\omega M/ F_{r-1}\omega M$ to $\sigma_{F_{r+1}M/ F_{r-1}M}\circ g\circ {\sigma_{Gr^F_{r+1}M}}^{-1}=\sigma_{F_{r+1}M/ F_{r-1}M}\circ g$ where as before, $\sigma_X$ means the action of $\sigma$ on $\omega X$ for an object $X$ of $\langle M\rangle^{\otimes}$ (note that $\sigma_{Gr^F_{r+1}M}=1$ because $\sigma\in\mathcal{U}(M,\omega)$). Writing endomorphisms of $F_{r+1}\omega M/F_{r-1}\omega M$ as $2\times 2$ matrices via our isomorphism $F_{r+1}\omega M/F_{r-1}\omega M \cong Gr^F_{r}\omega M \oplus Gr^F_{r+1}\omega M$, the matrix $\sigma_{F_{r+1}M/ F_{r-1}M}$ is obtained from $\sigma=(\sigma_{ij})$ by truncation to the part $r\leq i,j\leq r+1$, so that 
\[ \sigma_{F_{r+1}M/ F_{r-1}M}=\begin{pmatrix}
1 & f_{r,r+1}\\
0 & 1\end{pmatrix}.\]
Given any $g=(g_{r,r+1}, \lambda)\in \omega E_{r,r+1}$, we have
\[
\sigma\cdot g = \sigma_{F_{r+1}M/ F_{r-1}M}\circ g = \begin{pmatrix}
1 & f_{r,r+1}\\
0 & 1\end{pmatrix}\begin{pmatrix} g_{r,r+1} \\ \lambda\end{pmatrix}
\]
(the computation done in $Hom_\F(Gr^F_{r+1}\omega M, \displaystyle{\frac{F_{r+1}\omega M}{F_{r-1}\omega M}})$). Thus the morphism $\mathcal{U}(M,\omega)\rightarrow Hom_\F(Gr^F_{r+1}\omega M,Gr^F_{r}\omega M)$ of algebraic groups that corresponds to the action of $\mathcal{U}(M,\omega)$ on $E_{r,r+1}$ is simply given by $\sigma \mapsto f_{r,r+1}$. Taking logarithms, the image under $\omega$ of the morphism $\Psi_{V_{r,r+1}}(\mathcal{E}_r)$ is the map 
\[
\omega \fu(M)\rightarrow Hom_\F(Gr^F_{r+1}\omega M,Gr^F_{r}\omega M) \hspace{.5in} f\mapsto f_{r,r+1}.
\]
Upon recalling the description of $\pi_{r,r+1}$ from \S \ref{sec: defn of pi} this gives the result.
\medskip\par 
\item[(b)] By the functoriality of the map $\Psi_X$ in $X$ we have a commutative diagram
\[
\begin{tikzcd}[column sep = large]
Ext^1_{\langle M\rangle^{\otimes}}(\mathbbm{1}, V) \ar[d, equal] \arrow[r, "\Psi_V"] & Hom_\bT(\fu(M), V) \ar[d, equal] \\
\bigoplus\limits_r Ext^1_{\langle M\rangle^{\otimes}}(\mathbbm{1}, V_{r,r+1}) \arrow[r, "(\Psi_{V_{r,r+1}})"] & \bigoplus\limits_r Hom_\bT(\fu(M), V_{r,r+1})
\end{tikzcd}
\]
where the vertical identifications are via the canonical isomorphisms. Part (b) follows from part (a) since $\mathcal{E}$ (resp. $\pi$) corresponds to the tuple $(\mathcal{E}_{r,r+1})$ (resp. $(\pi_{r,r+1})$) under the vertical identifications.
\end{itemize}
\end{proof}

\subsection{Proof of Theorem \ref{thm: main thm}}\label{sec: proof of main thm}
We now deduce Theorem \ref{thm: main thm}. Recall from \eqref{eq2} in \S \ref{sec: filtration on u(M)} that $Gr^F_{-1}\fu(M)=\pi(\fu(M))$. Let $W$ be a subobject of $V$. We need to show that one has $\pi(\fu(M))\subset W$ if and only if the pushforward of $\mathcal{E}$ along the quotient map $V\rightarrow V/W$ lives in the subcategory $\langle Gr^FM\rangle^{\otimes}$. Tentatively, denote the quotient map $V\rightarrow V/W$ by $\varphi$. We have a commutative diagram
\begin{equation}\label{eq6}
\begin{tikzcd}[column sep = large]
Ext^1_{\langle M\rangle^{\otimes}}(\mathbbm{1}, V) \arrow[d, "\varphi_{\ast}"] \arrow[r, "\Psi_V"] & Hom_\bT(\fu(M), V) \arrow[d, "\varphi \, \circ -"] \\
Ext^1_{\langle M\rangle^{\otimes}}(\mathbbm{1}, V/W) \arrow[r, "\Psi_{V/W}"] &  Hom_\bT(\fu(M), V/W),
\end{tikzcd}
\end{equation}
where $\varphi_\ast$ denotes pushforward of extensions along $\varphi$. The kernel of the map $\Psi_{V/W}$ is the subgroup $Ext^1_{\langle Gr^FM\rangle^{\otimes}}(\mathbbm{1}, V/W)$ (see \S \ref{sec: Psi_X}). Thus $\varphi_\ast(\mathcal{E})$ is in the subcategory $\langle Gr^FM\rangle^{\otimes}$ if and only if $\Psi_{V/W}(\varphi_\ast(\mathcal{E}))=0$ if and only if $\varphi\circ (\Psi_V(\mathcal{E})) = 0$. Thanks to Proposition \ref{prop: main prop}, the latter is equivalent to $\pi(\fu(M))\subset W$. \hfill \qedsymbol{}

\begin{rem}
Suppose $Gr^FM$ is semisimple. Then Theorem \ref{thm: main thm} asserts that $Gr^F_{-1}\fu(M)$ is the smallest subobject of $V$ such that the pushforward of $\mathcal{E}$ along the quotient map $V\rightarrow V/Gr^F_{-1}\fu(M)$ splits. Since $Gr^FM$ is semisimple, this can be equivalently formulated as follows: $Gr^F_{-1}\fu(M)$ is the intersection of the kernels of all endomorphisms of $V$ that annihilate $\mathcal{E}$:
\[
Gr^F_{-1}\fu(M) \ = \bigcap\limits_{\stackrel{\phi\in End_\bT(V)}{\phi_\ast(\mathcal{E})=0}} \ker(\phi).
\]
\end{rem}


\subsection{A maximality criterion}\label{sec: maximality criteria}
In this subsection we use Theorem \ref{thm: main thm} to give a criterion for when $\fu(M)$ is equal to $F_{-1}\inEnd(M)$. Before we state the result, let us recall a definition due to Bertrand \cite{Ber01}: An extension $\mathcal{Z}$ of an object $Y$ by an object $X$ in $\bT$ is called totally nonsplit (or totally unsplit) if, considering $\mathcal{Z}$ as an extension of $\mathbbm{1}$ by $\inHom(Y,X)$, for every proper subobject $W$ of $\inHom(Y,X)$ the pushforward of $\mathcal{Z}$ along the quotient map $\inHom(Y,X)\rightarrow \inHom(Y,X)/W$ is nonsplit. If $\inHom(Y,X)$ is semisimple, then $\mathcal{Z}$ is totally nonsplit if and only if the annihilator of $\mathcal{Z}$ in $End_\bT(\inHom(Y,X))$ (for the $End_\bT(\inHom(Y,X))$-module $Ext^1_\bT(\mathbbm{1}, \inHom(Y,X))$) is trivial.

Theorem \ref{thm: main thm} has the following consequence:
\begin{prop}\label{prop: maximality criteria 1}
(a) The following statements are equivalent:
\begin{itemize}[wide]
\item[(i)] $\fu(M)=F_{-1}\inEnd(M)$
\item[(ii)] $Gr^F_{-1}\fu(M)=V$
\item[(iii)] If $W\subset V$ has the property that the pushforward of $\mathcal{E}$ along the quotient map $V\rightarrow V/W$ is in the subcategory $\langle Gr^FM\rangle ^{\otimes}$, then $W=V$. 
\end{itemize}
In particular, whether or not $\fu(M)=F_{-1}\inEnd(M)$ depends only on the extensions $\displaystyle{\frac{F_{r+1}M}{F_{r-1}M}}$ of $Gr^F_{r+1}M$ by $Gr^F_{r-1}M$.
\medskip\par 
\noindent (b) If $\fu(M)=F_{-1}\inEnd(M)$, then the extension $\mathcal{E}$ of $\mathbbm{1}$ by $V$ (and a fortiori, each $\mathcal{E}_r$ of $\mathbbm{1}$ by $V_{r,r+1}$) is totally nonsplit.
\medskip\par 
\noindent (c) If $Gr^FM$ is semisimple, then the statements (i)-(iii) of part (a) are also equivalent to the following statement:
\begin{itemize}[wide]
\item[(iv)] The extension $\mathcal{E}$ of $\mathbbm{1}$ by $V$ is totally nonsplit.
\end{itemize}
\end{prop}

\begin{proof}
(a) The equivalence (ii) $\Leftrightarrow$ (iii) is by Theorem \ref{thm: main thm}. The implication (i) $\Rightarrow$ (ii) is by diagram \eqref{eq2} in \S \ref{sec: filtration on u(M)}. In view of the same diagram and the fact that the derived algebra of $F_{-1}\inEnd(M)$ is $F_{-2}\inEnd(M)$, the implication (ii) $\Rightarrow$ (i) follows immediately from the following standard fact (see for instance, \cite[Lemma 7]{Ber13}): if $\fg$ is a nilpotent Lie algebra and $\fu$ is a Lie subalgebra of $\fg$ such that $\fu/(\fu\cap [\fg,\fg]) = \fg/[\fg,\fg]$, then $\fu=\fg$.
\medskip\par 
\noindent (b) This follows from the equivalent statement (iii): if the pushforward of $\mathcal{E}$ along $V\rightarrow V/W$ is split, then that pushforward is in $\langle Gr^FM\rangle ^{\otimes}$. 
\medskip\par 
\noindent (c) When $Gr^FM$ is semisimple, an extension of $\mathbbm{1}$ by a quotient of $V$ is in $\langle Gr^FM\rangle ^{\otimes}$ if and only if it splits. Thus in this case, (iii) of part (a) and statement (iv) are equivalent.
\end{proof}

The equivalence of statements (i) and (iii) in part (a) of Proposition \ref{prop: maximality criteria 1} will play an important role in \S \ref{sec: application to classification of motives with maximal uni rads}. The following definition (motivated by statement (iii) of Proposition \ref{prop: maximality criteria 1}(a)) will be convenient:

\begin{defn}\label{def: a totally disjoint extension}
Let $\mathbf{S}$ be a full Tannakian subcategory of $\bT$. Let $X$ be a nonzero object of $\bS$. We say an element $\mathcal{Z}\in Ext^1_\bT(\mathbbm{1},X)$ is {\it totally disjoint from $\mathbf{S}$} if for every proper subobject $W$ of $X$, the pushforward of $\mathcal{Z}$ along the quotient map $X\rightarrow X/W$ does not belong to $\mathbf{S}$. (Thus if $\mathbf{S}$ is semisimple, then $\mathcal{Z}$ is totally disjoint from $\mathbf{S}$ if and only if it is totally nonsplit.) 
\end{defn}

By Proposition \ref{prop: maximality criteria 1}(a), we have $\fu(M)=F_{-1}\inEnd(M)$ if and only if the extension $\mathcal{E}$ of $\mathbbm{1}$ by $V$ is totally disjoint from the subcategory $\langle Gr^FM\rangle^{\otimes}$.

\subsection{Variants - I}\label{sec: variants I}
In the remainder of \S \ref{sec: determination of pi(u(M))} we will discuss two variants of Theorem \ref{thm: main thm} and some consequences. This content will not be needed in \S \ref{sec: application to classification of motives with maximal uni rads}. A reader who wishes to skip to \S \ref{sec: application to classification of motives with maximal uni rads} may do so. 

One has the following variant (or generalization) of Theorem \ref{thm: main thm} about the projections of $\pi(\fu(M))$ to the direct summands of $V$:
\begin{prop}\label{prop: projections of pi(u(M)) to direct summands of V}
Let $I$ be a nonempty subset of $\{1,\ldots,k-1\}$. Set $V_I:=\bigoplus_{r\in I} V_{r,r+1}$ and $\pi_I:=(\pi_{r,r+1})_{r\in I}$, considered as a map from $\fu(M)$ to $V_I$. Let $\mathcal{E}_I:=(\mathcal{E}_r)_{r\in I}$, considered as an extension of $\mathbbm{1}$ by $V_I$. Then $\pi_I(\fu(M))$ is the smallest subobject of $V_I$ with the property that the pushforward of the extension $\mathcal{E}_I$ along the quotient map $V_I\rightarrow V_I/\pi_I(\fu(M))$ is an extension in the subcategory $\langle Gr^FM\rangle^{\otimes}$. In particular, if $Gr^FM$ is semisimple, then $\pi_I(\fu(M))$ is the smallest subobject of $V_I$ such that the pushforward of $\mathcal{E}_I$ to an extension of $\mathbbm{1}$ by $V_I/\pi_I(\fu(M))$ splits.
\end{prop}
\begin{proof}
By Proposition \ref{prop: main prop}(a) and functoriality of the maps $\Psi_X$ in $X$ we obtain $\Psi_{V_I}(\mathcal{E}_I) = \pi_I$. The proof of Proposition \ref{prop: projections of pi(u(M)) to direct summands of V} is now identical to the argument for Theorem \ref{thm: main thm} given in \S \ref{sec: proof of main thm}, with $V$, $\mathcal{E}$ and $\pi$ throughout being replaced by $V_I$, $\mathcal{E}_I$ and $\pi_I$, respectively.
\end{proof}

Consider a subquotient $F_nM/F_mM$ of $M$ with $0\leq m<n\leq k$. The filtration $F_\bullet$ induces a filtration on $F_nM/F_mM$. We can apply the constructions and results of the paper to this filtered object. We obtain an object $\fu(F_nM/F_mM)\subset F_{-1}\inEnd(F_nM/F_mM)$ whose image under any fiber functor is the Lie algebra of the subgroup of the Tannakian fundamental group of $F_nM/F_mM$ that acts trivially on $Gr^F(F_nM/F_mM)\cong \bigoplus_{m<r\leq n}Gr^F_r(M)$. We also have a natural map
\begin{equation}\label{eq33}
\fu(M) \rightarrow \fu(F_nM/F_mM)
\end{equation}
induced by the restriction map from the Tannakian group of $M$ to the Tannakian group of $F_nM/F_mM$. This map fits in a commutative diagram
\[
\begin{tikzcd}
\fu(M) \ar[r] \ar[d, phantom,"\bigcap"] & \fu(F_nM/F_mM) \ar[d, phantom,"\bigcap"] \\
F_{-1}\inEnd(M) \ar[r, twoheadrightarrow] &F_{-1}\inEnd(F_nM/F_mM),
\end{tikzcd}
\]
where the bottom arrow is the map that after applying a fiber functor $\omega$, sends a linear map $f\in End_\F(\omega M)$ in $F_{-1}\omega \inEnd(M)$ to the induced map in $End_\F(\omega F_nM/\omega F_mM)$ (this makes sense because $f$ is in $F_{-1}\omega \inEnd(M)$). One may ask if the top arrow is always surjective. As an application of Proposition \ref{prop: projections of pi(u(M)) to direct summands of V} we can see that this is not true in general, as we now explain.

For brevity, let $(m,n]$ denote the set of integers $\{r: m< r\leq n\}$. Consider the map $\pi$ for $F_nM/F_mM$ (see \S \ref{sec: defn of pi}). Including $F_nM/F_mM$ as a superscript in the notation to avoid confusion, this is a map
\[
\pi^{F_nM/F_mM} : F_{-1}\inEnd(F_nM/F_mM) \rightarrow V_{(m,n]} = \bigoplus\limits_{m<r<n} V_{r,r+1}
\]
(the notation $V_{(m,n]}$ is as in Proposition \ref{prop: projections of pi(u(M)) to direct summands of V}). There is a commutative diagram
\begin{equation}\label{eq14}
\begin{tikzcd}
\fu(M) \arrow[r, "\eqref{eq33}"] \arrow[rd, "\pi^M_{(m,n]}" '] & \fu(F_{n}M/F_{m}M) \arrow[d, "\pi^{F_{n}M/F_{m}M}"] \\
& V_{(m,n]},
\end{tikzcd}
\end{equation}
where $\pi^M$ means our $\pi$ from before (for $M$) and $\pi^M_{(m,n]}={(\pi_{r,r+1})}_{m<r<n}$ (again, following the notation of Proposition \ref{prop: projections of pi(u(M)) to direct summands of V}). By Proposition \ref{prop: projections of pi(u(M)) to direct summands of V}, $\pi^M_{(m,n]}(\fu(M))$ is the smallest subobject of $V_{(m,n]}$ such that the pushforward of $\mathcal{E}_{(m,n]}$ along the quotient $V_{(m,n]}\rightarrow V_{(m,n]}/\pi^M_{(m,n]}(\fu(M))$ belongs to $\langle Gr^FM\rangle^{\otimes}$. On the other hand, applying the same result (or Theorem \ref{thm: main thm}) to $F_{n}M/F_{m}M$ we see that $\pi^{F_{n}M/F_{m}M}(\fu(F_{n}M/F_{m}M))$ is the smallest subobject of $V_{(m,n]}$ such that the pushforward of $\mathcal{E}_{(m,n]}$ along the quotient $V_{(m,n]}\rightarrow V_{(m,n]}/\pi^{F_{n}M/F_{m}M}(\fu(F_{n}M/F_{m}M))$ belongs to $\langle Gr^F(F_nM/F_mM)\rangle^{\otimes}$. Already in the case $n-m=2$ with $n=r+1$ and $m=r-1$ (where $F_nM/F_mM$ has only two graded components $Gr^F_{r+1}M$ and $Gr^F_r M$), one can now easily use these characterizations to construct an example where the restriction map $\fu(M)\rightarrow \fu(F_{r+1}M/F_{r-1}M)$ is not surjective.

On the other hand, one has the following corollary:
\begin{cor}\label{cor: sufficient conditions for surjectivity of u(M) -> u(F_n M/F_m M)}
Suppose that $Gr^FM$ is semisimple.
\begin{itemize}[wide]
\item[(a)] Let $0\leq m<n\leq k$. Then with notation as above, we have
\[
\pi^M_{(m,n]}(\fu(M)) = \pi^{F_{n}M/F_{m}M}(\fu(F_{n}M/F_{m}M)).
\]
\item[(b)] For every $r$, the restriction map 
\[\fu(M)\rightarrow \fu(F_{r+1}M/F_{r-1}M)\] 
is surjective.
\item[(c)] Let $0\leq m<n\leq k$. If
\begin{equation}\label{eq27}
F_{-2}\fu(F_{n}M/F_{m}M) = [\fu(F_{n}M/F_{m}M),\fu(F_{n}M/F_{m}M)],
\end{equation}
then the restriction map
\[\fu(M)\rightarrow \fu(F_{n}M/F_{m}M)\] 
is surjective.
\end{itemize}
\end{cor}
\begin{proof}
\begin{itemize}[wide]
\item[(a)] Since $Gr^FM$ is semisimple, the conditions ``belongs to $\langle Gr^F M\rangle^{\otimes}$" and ``belongs to $\langle Gr^F(F_nM/F_m M)\rangle^{\otimes}$" for pushforwards of the extension $\mathcal{E}_{(m,n]}$ are both equivalent to splitting. Thus part (a) follows from Proposition \ref{prop: projections of pi(u(M)) to direct summands of V}.
\item[(b)] Setting $n=r+1$ and $m=r-1$, diagram \eqref{eq14} becomes
\[
\begin{tikzcd}
\fu(M) \arrow[r] \arrow[rd, "\pi_{r,r+1}" '] & \fu(F_{r+1}M/F_{r-1}M) \arrow[d, hookrightarrow,  "\pi^{F_{r+1}M/F_{r-1}M}"] \\
& V_{r,r+1} ,
\end{tikzcd}
\]
where $\pi^{\frac{F_{r+1}M}{F_{r-1}M}}$ is simply the natural inclusion of $\fu(F_{r+1}M/F_{r-1}M)$ in $F_{-1}\inEnd(\displaystyle{\frac{F_{r+1}M}{F_{r-1}M}})\cong V_{r,r+1}$. Thus the assertion follows from (a) and injectivity of $\pi^{F_{r+1}M/F_{r-1}M}$.
\item[(c)] The assumption \eqref{eq27} implies that $\pi^{F_{n}M/F_{m}M}(\fu(F_{n}M/F_{m}M))$ is the abelianization of the nilpotent Lie algebra $\fu(F_{n}M/F_{m}M)$. By part (a), the composition
\[
\fu(M) \rightarrow \fu(F_{n}M/F_{m}M) \rightarrow \fu(F_{n}M/F_{m}M)^{ab}
\]
is surjective. The claim follows from \cite[Lemma 7]{Ber13}.
\end{itemize}
\end{proof}

\begin{rem}
Corollary \ref{cor: sufficient conditions for surjectivity of u(M) -> u(F_n M/F_m M)}(b) is a special case of part (c) (as $F_{-2}\inEnd(F_{r+1}M/F_{r-1}M)$ is zero). I do not know whether or not for general $n,m$, the semisimplicity of $Gr^FM$ by itself (without condition \eqref{eq27}) is enough to guarantee the surjectivity of the restriction map $\fu(M)\rightarrow \fu(F_{n}M/F_{m}M)$.
\end{rem}

\begin{rem}
Suppose that after a relabelling of indices, the filtration $F_\bullet M$ is given by a functorial weight filtration $W_\bullet$ (see \S \ref{sec: application to classification of motives with maximal uni rads} for what this exactly means). Then the restriction maps $\fu(M)\rightarrow \fu(F_nM/F_mM)$ are all surjective, even if $Gr^FM$ is not semisimple. To see this, let $\omega_0$ be a fiber functor for $\bT$ and let $\omega=\omega_0\circ Gr^W$, where $Gr^W$ is the associated graded functor $\bT\rightarrow \bT$. Then for every object $N$ of $\bT$ the surjection $\mathcal{G}(N,\omega)\twoheadrightarrow \mathcal{G}(Gr^W N,\omega)$ admits a canonical section given by $Gr^W: \langle N\rangle^{\otimes}\rightarrow \langle Gr^W N\rangle^{\otimes}$. The surjectivity of $\fu(M)\rightarrow \fu(F_nM/F_mM)$ can be seen easily from this. 
\end{rem}

\subsection{Variants - II}\label{sec: variants II}
Suppose $Gr^FM$ is semisimple. In this case, we will give an equivalent formulation of Theorem \ref{thm: main thm} that will be useful to give a refinement of Proposition \ref{prop: maximality criteria 1}. Thanks to semisimplicity of $Gr^FM$, the object $\fu(F_{r+1}M/F_{r-1}M)$ is the smallest subobject of $V_{r,r+1}$ such that the pushforward of $\mathcal{E}_r$ along the quotient map $V_{r,r+1}\rightarrow V_{r,r+1}/\fu(F_{r+1}M/F_{r-1}M)$ splits (this is by the case $k=2$ of Theorem \ref{thm: main thm}, which was also proved earlier in \cite[Theorem 3.3.1]{EM1}). Equivalently, $\fu(F_{r+1}M/F_{r-1}M)$ is the smallest subobject of $V_{r,r+1}$ such that $\mathcal{E}_r$ is in the image of the map
\begin{equation}\label{eq15}
Ext^1_\bT(\mathbbm{1}, \fu(F_{r+1}M/F_{r-1}M)) \rightarrow Ext^1_\bT(\mathbbm{1}, V_{r,r+1})
\end{equation}
given by pushforward along the inclusion $\fu(F_{r+1}M/F_{r-1}M)\hookrightarrow V_{r,r+1}$. Since $V_{r,r+1}$ is semisimple, the map \eqref{eq15} is injective. With abuse of notation, we use the same notation for $\mathcal{E}_r$ and the element of $Ext^1_\bT(\mathbbm{1}, \fu(F_{r+1}M/F_{r-1}M))$ that pushes forward to it. 

Under the assumption of semisimplicity of $Gr^FM$, Theorem \ref{thm: main thm} can be equivalently formulated as follows:

\begin{prop}\label{prop: main thm in the ss case, variant version}
Suppose $Gr^FM$ is semisimple. Set
\[
\hat{V}:=\bigoplus\limits_r \, \fu(F_{r+1}M/F_{r-1}M) \subset V. 
\]
For each $r$, consider $\mathcal{E}_r$ as an extension of $\mathbbm{1}$ by $\fu(F_{r+1}M/F_{r-1}M)$, and consider $\mathcal{E}=(\mathcal{E}_r)$ as an extension of $\mathbbm{1}$ by $\hat{V}$. Then
$Gr^F_{-1}\fu(M)$ is the smallest subobject of $\hat{V}$ such that the pushforward of the extension $\mathcal{E}$ along the quotient map $\hat{V}\rightarrow \hat{V}/Gr^F_{-1}\fu(M)$ splits.
\end{prop}

As an immediate consequence of Proposition \ref{prop: main thm in the ss case, variant version} one gets the following result, which refines Proposition \ref{prop: maximality criteria 1} in the case when $Gr^FM$ is semisimple and $\hat{V}\neq V$.
\begin{cor}
Suppose $Gr^FM$ is semisimple. Then the following statements are equivalent:
\begin{itemize}
\item[(i)] $Gr^F_{-1}\fu(M)=\hat{V}$
\item[(ii)] The extension $\mathcal{E}$, considered as an extension of $\mathbbm{1}$ by $\hat{V}$, is totally nonsplit (see \S \ref{sec: maximality criteria} to recall what this means).
\end{itemize}
\end{cor}

Note that for semisimple $Gr^FM$, one has $\hat{V}=V$ if and only if for every $r$, the extension $\mathcal{E}_r\in Ext^1_\bT(Gr^F_{r+1}M, Gr^F_{r}M)$ is totally nonsplit.

\section{Classification of objects with maximal $\fu$ in the case of a weight filtration}\label{sec: application to classification of motives with maximal uni rads}

\subsection{Statement of the result}
In the remainder of the paper, we assume that our Tannakian category $\bT$ is equipped with a functorial filtration $W_\bullet$ (called the weight filtration) with similar properties to the weight filtration on rational mixed Hodge structures: indexed by $\ZZ$, increasing and finite on every object, functorial, exact, and compatible with the tensor structure. More explicitly, this means that every object $M$ of $\bT$ is equipped with a filtration $(W_n M)_{n\in \ZZ}$ with $W_{n-1}M\subset W_nM$ for all $n\in\ZZ$, $W_n M = 0$ for $n\ll 0$, and $W_n M = M$ for $n\gg 0$, such that the following conditions (i) - (iii) hold: (i) For every objects $M$ and $N$ and morphism $f:M\rightarrow N$ in $\bT$, we have $f(W_n M)\subset W_n N$ for every $n$; (ii) the functors $W_n: \bT\rightarrow \bT$ (defined by $M\mapsto W_n M$ on objects and acting on a morphism $f: M\rightarrow N$ by restricting it to $f: W_nM\rightarrow W_nN$) are exact for every $n$; and (iii) for every objects $M$ and $N$ and every integer $n$,
\[
W_n(M\otimes N) = \sum\limits_{p+q=n}W_pM\otimes W_qN.
\]
We call an object $M$ of $\bT$ pure if there exists an integer $n$ such that $W_{n-1}M=0$ and $W_nM=M$. A weight of an object $M$ is an integer $n$ such that $W_{n-1}M\neq W_nM$. A nonzero pure object is an object with exactly one weight.

For every nonzero object $M$ of $\bT$, we apply the constructions and results of the paper by taking the increasing filtration $(F_r M)_{r\in\ZZ}$ to be given by the weight filtration: more precisely, if the weights of $M$ are $p_1<\cdots<p_k$, set $F_rM:=W_{p_r}M$ for $1\leq r\leq k$ and $F_0M=0$. We thus obtain an object 
\[\fu(M)\subset F_{-1}\inEnd(M)=W_{-1}\inEnd(M).\]

As an application of Theorem \ref{thm: main thm} we will prove the following result about the structure of the set of isomorphism classes of objects $M$ of $\bT$ with a given associated graded and maximal $\fu(M)$:

\begin{thm}\label{thm: isomorphism classes of objects with given Gr and maximal u}
Fix objects $A_1,\ldots, A_k$ of $\bT$, where $A_r$ is nonzero and pure of weight $p_r$, and $p_1<\cdots < p_k$. Set $A := \bigoplus_r A_r$. Let $S^\ast(A)$ be the set of equivalence classes of objects $M$ of $\bT$ such that $Gr^WM$ is isomorphic to $A$ and $\fu(M)=W_{-1}\inEnd(M)$, where two such objects $M$ and $M'$ are considered equivalent if they are isomorphic in $\bT$ (note that we do not keep any trace of a choice of isomorphism between $Gr^WM$ and $A$). Then there exist sets $S^\ast_\ell(A)$ ($1\leq \ell\leq k-1$) and maps
\begin{equation}\label{eq20}
S^\ast(A)\cong S^\ast_{k-1}(A) \rightarrow S^\ast_{k-2}(A) \rightarrow\cdots \rightarrow S^\ast_1(A)
\end{equation}
such that the following statements hold:
\begin{itemize}[wide]
\item[(i)] For $2\leq \ell\leq k-1$, every nonemply fiber of the map $S_\ell^\ast(A)\rightarrow S_{\ell-1}^\ast(A)$ is a torsor over
\[
\bigoplus\limits_{r} Ext^1_\bT(A_{r+\ell},A_r)
\]
(where $r$ runs from $1$ to $k-\ell$).
\item[(ii)] Let $2\leq \ell\leq k-1$. Then the map $S_\ell^\ast(A)\rightarrow S_{\ell-1}^\ast(A)$ is surjective if
\[
\bigoplus\limits_{r} Ext^2_\bT(A_{r+\ell},A_r)
\]
vanishes.
\item[(iii)] (Description of $S^\ast_1(A)$) Let
\[
\bigm(\bigoplus_r Ext^1_\bT(A_{r+1},A_r)\bigm)^\ast
\]
be the subset of $\bigoplus_r Ext^1_\bT(A_{r+1},A_r)$ consisting of every extension tuple $\mathcal{E}=(\mathcal{E}_r)$ such that considering $\mathcal{E}$ as an element of $Ext^1_\bT(\mathbbm{1}, \bigoplus_r \inHom(A_{r+1},A_r) )$ in the natural way, $\mathcal{E}$ is totally disjoint from the subcategory $\langle A\rangle^{\otimes}$ (see Definition \ref{def: a totally disjoint extension}). Then there is a canonical bijection 
\[S^\ast_1(A) \ \cong \ \bigm(\bigoplus_r Ext^1_\bT(A_{r+1},A_r)\bigm)^\ast \bigm/ Aut(A),\]
where the action of $Aut(A)=\prod_r Aut(A_r)$ on extensions is by pushforwards and pullbacks (i.e., $\sigma=(\sigma_r)\in \prod_r Aut(A_r)$ sends $(\mathcal{E}_r)$ to $((\sigma_{r+1}^{-1})^\ast (\sigma_{r})_\ast \mathcal{E}_r)$).
\end{itemize}
\end{thm}

This result was proved in \cite[\S 5, Theorem 5.4.3]{Es1} in the special case when $A$ is ``graded-independent" (in the sense of Definition 5.3.1 of \cite{Es1}) and semisimple. Our goal in the rest of the paper is to establish Theorem \ref{thm: isomorphism classes of objects with given Gr and maximal u}. This will be done in \S \ref{sec: proof of classification thm} after we review some relevant constructions of \cite{Es1} in \S \ref{sec: recollections on gen exts} below. In what follows until the end of the paper, $A$ (and $k$, the $A_r$ and $p_r$) are fixed and as described in the statement of Theorem \ref{thm: isomorphism classes of objects with given Gr and maximal u}.

\subsection{Recollections on generalized extensions}\label{sec: recollections on gen exts}
Following \cite{Es1}, we denote by $S'(A)$ the set of equivalence classes of pairs $(M,\phi)$ consisting of an object $M$ of $\bT$ and an isomorphism $\phi: Gr^WM\rightarrow A$, where two such pairs $(M,\phi)$ and $(M',\phi')$ are considered equivalent if there exists a morphism (automatically, an isomorphism) $f: M\rightarrow M'$ such that $\phi'\circ Gr^Wf=\phi$. Denote by $S(A)$ the set of isomorphism classes of objects $M$ of $\bT$ such that $Gr^WM$ is isomorphic to $A$. Note that no trace of a choice of an isomorphism $Gr^WM\rightarrow A$ is kept in $S(A)$ (in contrast to the set $S'(A)$). There is a natural action of $Aut(A)$ on $S'(A)$ (given by twisting $\phi$ in a pair $(M,Gr^WM\xrightarrow{\phi} A)$), and the set $S(A)$ can be identified with the quotient of $S'(A)$ by $Aut(A)$. The set $S^\ast(A)$ of Theorem \ref{thm: isomorphism classes of objects with given Gr and maximal u} is contained in $S(A)$.

We introduced a new approach to study $S'(A)$ and $S(A)$ in \cite[\S 3]{Es1}. This approach, which in loc. cit. we called induction on the {\it level}, appears to have better naturalness properties than the obvious approach of induction on the number of weights $k$ (see \cite[Remark 3.3.2]{Es1}). The inductive approach on the level is a natural generalization of the formalism of blended extensions ({\it extensions panach\'{e}es} \cite{Gr68}). Here we will briefly sketch this approach to the extent needed to prove Theorem \ref{thm: isomorphism classes of objects with given Gr and maximal u}. We refer the reader to \S 3 and \S 4 of \cite{Es1} for more details and the proofs of everything we will say below in this subsection.

We constructed sets $S'_\ell(A)$ for $1\leq \ell\leq k-1$ and maps
\begin{equation}\label{eq16}
S'(A)\cong S'_{k-1}(A) \rightarrow S'_{k-2}(A) \rightarrow\cdots \rightarrow S'_1(A)\cong \bigoplus\limits_{r}Ext^1_\bT(A_{r+1},A_r)
\end{equation}
such that for every $\ell$, every nonempty fiber of $S'_\ell(A)\rightarrow S'_{\ell-1}(A)$ is canonically a torsor over
\[
\bigoplus\limits_{r} Ext^1_\bT(A_{r+\ell},A_r).
\]
Moreover, the map $S'_\ell(A)\rightarrow S'_{\ell-1}(A)$ is surjective if
\[
\bigoplus\limits_{r} Ext^2_\bT(A_{r+\ell},A_r)
\]
vanishes. We also constructed quotients $S_\ell(A)$ of $S'_\ell(A)$ for $1\leq \ell\leq k-1$ such that \eqref{eq16} descends to a sequence of maps
\begin{equation}\label{eq17}
S(A)\cong S_{k-1}(A) \rightarrow S_{k-2}(A) \rightarrow\cdots \rightarrow S_1(A)\cong \  \bigm(\bigoplus\limits_{r} Ext^1_\bT(A_{r+1},A_r)\bigm)\bigm/Aut(A), 
\end{equation}
where the action of the automorphism group $Aut(A)=\prod_r Aut(A_r)$ on a tuple of extension classes is by pushforwards and pullbacks, as described at the end of Theorem \ref{thm: isomorphism classes of objects with given Gr and maximal u}(iii). (See Theorem 3.3.1(a-d) of \cite{Es1} for all of this.)

We will sketch the constructions of \eqref{eq16} and \eqref{eq17} as they will be relevant. We shall use the word ``depth" for what was called ``level" in \cite{Es1}, since the former might be a better choice of wording. Given an object $M$ and an isomorphism $\phi: Gr^WM\rightarrow A$, setting $M_{m,n}:=W_{p_n}M/W_{p_m}M$ for $0\leq m<n\leq k$ with $p_0:=p_1-1$ (so that $W_{p_0}M=0$), the natural inclusions and projections between the $M_{m,n}$ give rise to a diagram
\begin{equation}\label{eq18}
\begin{tikzcd}[column sep=small, row sep=small]
A_1         \arrow[r, hookrightarrow] & M_{0,2} \arrow[r, hookrightarrow] \arrow[d, twoheadrightarrow] & M_{0,3} \arrow[r, hookrightarrow] \arrow[d, twoheadrightarrow] & M_{0,4} \arrow[r, hookrightarrow] \arrow[d, twoheadrightarrow] & \cdots \arrow[r, hookrightarrow] & M_{0,k-1} \arrow[r, hookrightarrow] \arrow[d, twoheadrightarrow] & M_{0,k} \arrow[d, twoheadrightarrow] \\
 & A_2 \arrow[r, hookrightarrow] & M_{1,3} \arrow[d, twoheadrightarrow] \arrow[r, hookrightarrow] & M_{1,4} \arrow[d, twoheadrightarrow] \arrow[r, hookrightarrow] & \cdots \arrow[r, hookrightarrow] & M_{1,k-1} \arrow[r, hookrightarrow] \arrow[d, twoheadrightarrow] & M_{1,k} \arrow[d, twoheadrightarrow] \\
 &  & A_3 \arrow[r, hookrightarrow] & M_{2,4} \arrow[r, hookrightarrow] \arrow[d, twoheadrightarrow] & \cdots & \vdots & \vdots \\
& &  & A_4 \arrow[r, hookrightarrow] & \cdots & & \\
 & & & & \ddots & \arrow[d, twoheadrightarrow] & \arrow[d, twoheadrightarrow] \\
&&&&& A_{k-1} \arrow[r,hookrightarrow] & M_{k-2,k} \arrow[d, twoheadrightarrow] \\
&&&&& & A_{k},
\end{tikzcd}
\end{equation}
where we have used the isomorphism $\phi$ to replace $M_{r-1,r}$ by $A_r$ for every $r$. A generalized extension of depth $k-1$ (or {\it full depth}) of $A$ is an abstract version of the data of the diagram above (see Definition 3.4.1 of \cite{Es1}). More generally, for $1\leq \ell\leq k-1$, a generalized extension of depth $\ell$ of $A$ is an abstract version of a truncated version of the diagram above, consisting of the first $\ell$ diagonals above the diagonal formed by the $A_r$. (By the $i$-th diagonal here we mean the objects $M_{m,n}$ with $n-m=i+1$.) The precise definition is recalled below. Note that in this definition and in what follows after, for convenience, in the context of indices we will often use the phrase ``in the eligible range" to mean in the range in which the indices in question makes sense.

\begin{defn}[Definition 3.4.2 of \cite{Es1}]\label{def: gen ext level l} Let $1\leq \ell\leq k-1$. A generalized extension of depth $\ell$ of $A$ consists of the data of an object $M_{m,n}$ of $\bT$ for each pair $(m,n)$ of integers with $0\leq m<n\leq k$ and $n-m\leq \ell+1$, with $M_{r-1,r}=A_r$ for all $1\leq r\leq k$, together with the data of a surjective morphism $M_{m,n}\rightarrow M_{m+1,n}$ and an injective morphism $M_{m,n-1}\rightarrow M_{m,n}$ for every $m$ and $n$ in the eligible ranges, such that the following two axioms hold:
\begin{itemize}[wide]
\item[(i)] Every diagram of the form
\[
\begin{tikzcd}
M_{m,n-1} \arrow[r, hookrightarrow] \arrow[d, twoheadrightarrow] & M_{m,n} \arrow[d, twoheadrightarrow] \\
M_{m+1,n-1} \arrow[r, hookrightarrow] & M_{m+1,n} 
\end{tikzcd}
\]
(with the maps as in the given data) commutes.
\item[(ii)] Every diagram of the form
\begin{equation}\label{eq19}
\begin{tikzcd}
   0 \arrow[r] &  M_{m,n-1}\arrow[r, ] & M_{m,n} \arrow[r, ] &  A_n  \arrow[r] & 0
\end{tikzcd}
\end{equation}
is an exact sequence. Here, the morphism $M_{m,n}\rightarrow A_n$ is the composition
\[
M_{m,n} \twoheadrightarrow M_{m+1,n} \twoheadrightarrow M_{m+2,n} \twoheadrightarrow \cdots \twoheadrightarrow M_{n-1,n}=A_n.  
\]
\end{itemize}
\end{defn}
\medskip\par 

We visualize a generalized extension of depth $\ell=k-1$ (resp. $\ell<k-1$) by a diagram of the form \eqref{eq18} (resp. a truncated version of the diagram that includes the $\ell$ diagonals above the $A_r$). For example, the diagram below shows a generalized extension of depth $\ell=2$ in the case $k=4$:
\[
\begin{tikzcd}[column sep=small, row sep=small]
A_1         \arrow[r, hookrightarrow] & M_{0,2} \arrow[r, hookrightarrow] \arrow[d, twoheadrightarrow] & M_{0,3} \arrow[d, twoheadrightarrow] &  \\
 & A_2 \arrow[r, hookrightarrow] & M_{1,3} \arrow[d, twoheadrightarrow] \arrow[r, hookrightarrow] & M_{1,4} \arrow[d, twoheadrightarrow]  \\
 &  & A_3 \arrow[r, hookrightarrow] & M_{2,4} \arrow[d, twoheadrightarrow] \\
& &  & A_4.
\end{tikzcd}
\]

Given two generalized extensions $(M_\db)$ and $(N_\db)$ of the same depth, a morphism between them is defined as a collection of morphisms $M_{m,n}\rightarrow N_{m,n}$ in $\bT$ that commute with the structure arrows of $(M_\db)$ and $(N_\db)$. For each $1\leq \ell\leq k-1$, the set $S_\ell(A)$ of \eqref{eq17} (resp. $S'_\ell(A)$ of \eqref{eq16}) is defined as the collection of all generalized extensions of depth $\ell$ of $A$ modulo the equivalence relation given by declaring two generalized extensions to be equivalent if there exists an isomorphism between them (resp. an isomorphism between them that is identity on the $A_r$). In particular, there is a natural surjection $S'_\ell(A)\rightarrow S_\ell(A)$ for each $\ell$. The identification $S'_1(A) \cong \bigoplus_r Ext^1_\bT(A_{r+1},A_r)$ sends the equivalence class of a generalized extension $(M_{m,n})_{n-m\leq 2}$ of depth 1 (i.e., with just one diagonal above the $A_r$) to the element $(\mathcal{E}_r)$ of $\bigoplus_r Ext^1_\bT(A_{r+1},A_r)$, where $\mathcal{E}_r$ is given by
\begin{equation}\label{eq34}
\begin{tikzcd}
0 \ar[r] & A_r \ar[r] & M_{r-1,r+1} \ar[r] & A_{r+1} \ar[r] & 0
\end{tikzcd}
\end{equation}
(the maps being structure arrows of our generalized extension; note that this is indeed an extension by axiom (ii) of Definition \ref{def: gen ext level l}). The canonical bijection $S'(A)\cong S_{k-1}'(A)$ is given by sending the equivalence class of a pair $(M,Gr^WM\xrightarrow{\phi} A)$ to the class of the generalized extension \eqref{eq18} of full depth with $M_{m,n}:=W_{p_n}M/W_{p_m}M$ and $Gr^W_{p_r}M$ having been replaced by $A_r$ via $\phi$. This descends to a canonical bijection $S(A)\cong S_{k-1}(A)$. The maps $S'_\ell(A)\rightarrow S'_{\ell-1}(A)$ and $S_\ell(A)\rightarrow S_{\ell-1}(A)$ are given by {\it truncation}, i.e., by erasing the top diagonal ( = the $\ell$-th diagonal above the $A_r$) of a generalized extension of depth $\ell$. With abuse of notation, we denote the truncation maps $S'_\ell(A)\rightarrow S'_{\ell-1}(A)$ and $S_\ell(A)\rightarrow S_{\ell-1}(A)$ both by $\Theta_\ell$.

Notably, the group actions that describe how the fibers of the truncation map $S_\ell(A)\rightarrow S_{\ell-1}(A)$ descend from those of $S'_\ell(A)\rightarrow S'_{\ell-1}(A)$ were described in \cite{Es1} (see \S 3.10 therein). In particular, by Proposition 4.2.3 of \cite{Es1} (also see Remark 4.2.5 therein), given $\epsilon\in S_{\ell-1}(A)$ and any $\epsilon'\in S'_{\ell-1}(A)$ above $\epsilon$, if $\epsilon$ (equivalently, $\epsilon'$) is the class of a generalized extension $(M_\db)$ of depth $\ell-1$ {\it such that all the extensions \eqref{eq19} are totally nonsplit}, then the fiber of $S'_\ell(A)\rightarrow S'_{\ell-1}(A)$ above $\epsilon'$ descends bijectively to the fiber of $S_\ell(A)\rightarrow S_{\ell-1}(A)$ above $\epsilon$:
\begin{equation}\label{eq21}
\begin{tikzcd}
\Theta_\ell^{-1}(\epsilon') \ar[d, equal] \arrow[r, hookrightarrow] & S'_\ell(A) \arrow[d, twoheadrightarrow]   \ar[r, "\Theta_\ell"] & S'_{\ell-1}(A) \arrow[d, twoheadrightarrow] \\
\Theta_\ell^{-1}(\epsilon) \arrow[r, hookrightarrow] & S_\ell(A) \arrow[r, "\Theta_\ell"] & S_{\ell-1}(A).
\end{tikzcd}
\end{equation}
This last fact will be important for us in \S \ref{sec: proof of classification thm}.

As mentioned earlier, the reader can consult \cite[\S 3 and \S 4]{Es1} for the details and the proofs of all of what was discussed above. Note that the diagrams of generalized extensions above are the reflected versions of those in \cite{Es1} (reflected over the diagonal of the $A_r$), with the role of horizontal and vertical arrows switched. The pictures as shown here are more consistent with our matrix notation for elements of Tannakian groups of objects with increasing filtrations (where the matrices are upper triangular).

\subsection{Proof of Theorem \ref{thm: isomorphism classes of objects with given Gr and maximal u}}\label{sec: proof of classification thm}
We can now prove Theorem \ref{thm: isomorphism classes of objects with given Gr and maximal u}. The notation ($A_r$, $p_r$, $A$, $S^\ast(A)$ and $(\bigoplus_r Ext^1_\bT(A_{r+1},A_r))^\ast$) is as introduced in the statement of the theorem. We shall set
\[
V:=\bigoplus_{r} \inHom(A_{r+1},A_r).
\]
This notation is consistent with the notation of \S \ref{sec: determination of pi(u(M))} in the following way: If $M$ is an object of $\bT$ with $Gr^WM$ isomorphic to $A$, a choice of an isomorphism $Gr^WM\rightarrow A$ allows us to identify $V$ defined here with its namesake in \S \ref{sec: determination of pi(u(M))} and Theorem \ref{thm: main thm} for the filtration $(F_r M)$ on $M$ given by  $F_rM=W_{p_r}M$.

We start the proof by a lemma.

\begin{lemma}\label{lem: totally disjoint invariant under the action of Aut(A)}
\begin{itemize}[wide]
\item[(a)]
Let $\phi: A'\rightarrow A$ be an isomorphism in $\bT$, where $A'=\bigoplus_{1\leq r\leq k} A'_r$ and for every $r$, the component $A'_r$ is nonzero pure of weight $p_r$. Set $V':= \bigoplus_r\inHom(A'_{r+1},A'_r)$. Define 
\[\bigm(\bigoplus_r  Ext^1_\bT(A'_{r+1},A'_r)\bigm)^\ast\] 
in the same way as we defined $(\bigoplus_r Ext^1_\bT(A_{r+1},A_r))^\ast$ in statement (iii) of Theorem \ref{thm: isomorphism classes of objects with given Gr and maximal u}, i.e., $(\bigoplus_r Ext^1_\bT(A'_{r+1},A'_r))^\ast$ is the subset of $\bigoplus_r  Ext^1_\bT(A'_{r+1},A'_r)$ consisting of every tuple $\mathcal{E}'=(\mathcal{E}'_r)$ with the property that, considered as an extension of $\mathbbm{1}$ by $V'$ in the natural way, $\mathcal{E}'$ is totally disjoint from $\langle A\rangle^{\otimes}$ (in the sense of Definition \ref{def: a totally disjoint extension}). Let
\begin{equation}\label{eq30}
\bigoplus_r  Ext^1_\bT(A'_{r+1},A'_r) \rightarrow \bigoplus_r  Ext^1_\bT(A_{r+1},A_r)
\end{equation}
be the isomorphism given by pushforwards and pullbacks of extensions along $\phi$ (more explicitly, if $\phi=(\phi_r)$ with $\phi_r:A'_r\xrightarrow{\simeq} A_r$ and $\mathcal{E}'=(\mathcal{E}'_r)\in \bigoplus_r  Ext^1_\bT(A'_{r+1},A'_r)$, then \eqref{eq30} sends $\mathcal{E}'$ to $((\phi_{r+1}^{-1})^\ast(\phi_r)_\ast(\mathcal{E}'_r)$). Then the map \eqref{eq30} sends $(\bigoplus_r Ext^1_\bT(A'_{r+1},A'_r))^\ast$ bijectively to $(\bigoplus_r Ext^1_\bT(A_{r+1},A_r))^\ast$. 
\medskip\par 
\item[(b)] The subset $(\bigoplus_r Ext^1_\bT(A_{r+1},A_r))^\ast$ of $\bigoplus_r Ext^1_\bT(A_{r+1},A_r)$ is closed under the action of $Aut(A)$ on the latter space by pushforwards and pullbacks.
\end{itemize}
\end{lemma}
\begin{proof}
Part (b) is immediate from (a) so we only need to prove the latter. Thanks to the functoriality properties of $\inHom$, the isomorphism $\phi=(\phi_r): A'\rightarrow A$ induces an isomorphism $V'\rightarrow V$, which we denote also by $\phi$; given a fiber functor $\omega$, the isomorphism $\omega\phi: \omega V'\rightarrow \omega V$ sends $f\in Hom_\F(\omega A'_{r+1},\omega A'_r)$ to $(\omega\phi_r)\circ f\circ (\omega\phi_{r+1})^{-1}$. We have a commutative diagram
\[
\begin{tikzcd}
\bigoplus\limits_r Ext^1_\bT(A'_{r+1},A'_r) \ar[r] \ar[d, "\eqref{eq30}"] & Ext^1_\bT(\mathbbm{1}, V') \ar[d, "\phi_\ast"]\\
\bigoplus\limits_r Ext^1_\bT(A_{r+1},A_r) \ar[r] & Ext^1_\bT(\mathbbm{1}, V),
\end{tikzcd}
\]
where the horizontal arrows are the canonical isomorphisms and $\phi_\ast$ is pushforward along $\phi: V'\rightarrow V$. Thus to establish the lemma it suffices to verify that $\mathcal{E}'\in Ext^1_\bT(\mathbbm{1}, V')$ is totally disjoint from $\langle A\rangle^{\otimes}$ if and only if $\phi_\ast(\mathcal{E}')\in Ext^1_\bT(\mathbbm{1}, V)$ is such. To verify this, let $\mathcal{E}'$ be represented by an extension
\[
\begin{tikzcd}
0 \ar[r] & V' \ar[r, "\iota"] & E' \ar[r] & \mathbbm{1} \ar[r] & 0.
\end{tikzcd}
\]
Given $W\subset V$, consider the following two extension classes: (i) the pushforward of $\mathcal{E}'$ along the quotient $V'\rightarrow V'/\phi^{-1}(W)$, and (ii) the pushforward of $\phi_\ast(\mathcal{E}')$ along the quotient $V\rightarrow V/W$. These two extension classes are respectively given by the first and second rows of the commutative diagram
\[
\begin{tikzcd}
0 \ar[r] & V'/\phi^{-1}(W) \ar[r, "\iota"] \ar[d, "\phi"] & E'/\iota\phi^{-1}(W) \ar[r] \ar[d, equal] & \mathbbm{1} \ar[r] \ar[d, equal] & 0 \\
0 \ar[r] & V/W \ar[r, "\iota\phi^{-1}"] & E'/\iota\phi^{-1}(W) \ar[r] & \mathbbm{1} \ar[r] & 0. 
\end{tikzcd}
\]
In particular, the two extension classes have isomorphic middle objects and hence are in $\langle A\rangle^{\otimes}$ at the same time. Thus $\mathcal{E}'$ is totally disjoint from $\langle A\rangle^{\otimes}$ if and only if $\phi_\ast(\mathcal{E}')$ is such.
\end{proof}

With Lemma \ref{lem: totally disjoint invariant under the action of Aut(A)} in hand, we proceed to give the definition of the sets $S^\ast_\ell(A)$ of Theorem \ref{thm: isomorphism classes of objects with given Gr and maximal u}. We shall use the notation introduced in \S \ref{sec: recollections on gen exts} for generalized extensions. For $1\leq \ell\leq k-1$, the set $S^\ast_\ell(A)$ is the pre-image of
\begin{equation}\label{eq29}
\bigm(\bigoplus_r Ext^1_\bT(A_{r+1},A_r)\bigm)^\ast \bigm/ Aut(A) \ \ \subset \ \ \bigm(\bigoplus_r Ext^1_\bT(A_{r+1},A_r)\bigm) \bigm/ Aut(A)
\end{equation}
(which makes sense by Lemma \ref{lem: totally disjoint invariant under the action of Aut(A)}(b)) under the composition
\[
S_\ell(A)\rightarrow \cdots \rightarrow S_1(A) \stackrel{(\dagger)}{\cong} \bigm(\bigoplus_r Ext^1_\bT(A_{r+1},A_r)\bigm) \bigm/ Aut(A)
\]
in \eqref{eq17}. The requirement given in statement (iii) of Theorem \ref{thm: isomorphism classes of objects with given Gr and maximal u} is met thanks to the bijection $(\dagger)$. The map $S^\ast_\ell(A)\rightarrow S^\ast_{\ell-1}(A)$ in the statement of the theorem is just the restriction of the truncation map $S_\ell(A)\rightarrow S_{\ell-1}(A)$ in \eqref{eq17}.
\medskip\par 
Our next task is to argue that the bijection $S(A)\xrightarrow{\simeq} S_{k-1}(A)$ in \eqref{eq17} restricts to a bijection between $S^\ast(A)$ and $S^\ast_{k-1}(A)$. By construction, the composition 
\[S(A)\rightarrow S_{k-1}(A)\rightarrow \cdots \rightarrow S_1(A) \cong \bigm(\bigoplus\limits_{r} Ext^1_\bT(A_{r+1},A_r)\bigm)\bigm/Aut(A)\]
is the map that sends the isomorphism class of an object $M$ to the $Aut(A)$-orbit of the tuple of extension classes
\[
(W_{p_{r+1}}M/W_{p_{r-1}} M) \in \bigoplus_r Ext^1_\bT(Gr^W_{p_{r+1}}M, Gr^W_{p_{r}}M),
\]
with each $Gr^W_{p_r}M$ being replaced by $A_r$ via a choice of isomorphism $\phi: Gr^WM\rightarrow A$. By Proposition \ref{prop: maximality criteria 1}(a) and Lemma \ref{lem: totally disjoint invariant under the action of Aut(A)}(a), this sends $S^\ast(A)$ into \eqref{eq29}. Thus the bijection $S(A)\rightarrow S_{k-1}(A)$ restricts to an injection $S^\ast(A)\rightarrow S^\ast_{k-1}(A)$. On the other hand, given $\epsilon\in S^\ast_{k-1}(A)$, let $M$ be an object of $\bT$ with $Gr^WM\simeq A$ such that the class of $M$ in $S(A)$ is sent to $\epsilon$ by the bijection $S(A)\rightarrow S_{k-1}(A)$. Since the truncation of $\epsilon$ to depth 1 lives in \eqref{eq29}, the element $(W_{p_{r+1}}M/W_{p_{r-1}} M)$ of
\begin{equation}\label{eq31}
Ext^1_\bT\bigm(\mathbbm{1}, \bigoplus_r \inHom(Gr^W_{p_{r+1}}M, Gr^W_{p_{r}}M)\bigm)
\end{equation}
is totally disjoint from $\langle A\rangle^{\otimes}$ (first, we see this for the tuple $(W_{p_{r+1}}M/W_{p_{r-1}} M)$ considered as an element of $Ext^1_\bT(\mathbbm{1}, V)$ via a choice of isomorphism $\phi: Gr^WM\rightarrow A$, and then by Lemma \ref{lem: totally disjoint invariant under the action of Aut(A)}(a) for $(W_{p_{r+1}}M/W_{p_{r-1}} M)$ as an element of \eqref{eq31}). It now follows from Proposition \ref{prop: maximality criteria 1}(a) that $\fu(M)=W_{-1}\inEnd(M)$, so that the class of $M$ is indeed in $S^\ast(A)$. We have proved that the bijection $S(A)\cong S_{k-1}(A)$ restricts to a bijection $S^\ast(A)\cong S^\ast_{k-1}(A)$.
\medskip\par 
It remains to verify requirements (i) and (ii) of the statement of the theorem. Let $\epsilon\in S^\ast_{\ell-1}(A)$. By construction, the fiber of $S^\ast_\ell(A)\rightarrow S^\ast_{\ell-1}(A)$ above $\epsilon$ is the same as the fiber of $S_\ell(A)\rightarrow S_{\ell-1}(A)$ above $\epsilon$. We will see that if $(M_\db)$ is a generalized extension of depth $\ell$ representing $\epsilon$, then all the extensions \eqref{eq19} are totally nonsplit. Taking this for granted for the moment, choosing $\epsilon'\in S'_{\ell-1}(A)$ above $\epsilon$, by Proposition 4.2.3 of \cite{Es1} (also see Remark 4.2.5 therein) the fiber of $S_\ell(A)\rightarrow S_{\ell-1}(A)$ above $\epsilon$ is in a canonical bijection with the fiber of $S'_\ell(A)\rightarrow S'_{\ell-1}(A)$ above $\epsilon'$ (see diagram \eqref{eq21}). Statements (i) and (ii) of Theorem \ref{thm: isomorphism classes of objects with given Gr and maximal u} thus follow from the structure of the fibers of the truncation map $S'_\ell(A)\rightarrow S'_{\ell-1}(A)$ given in \cite[Theorem 3.3.1 (b,c)]{Es1} and recalled earlier in \S \ref{sec: recollections on gen exts}.
\medskip\par 
To complete the proof, let $(M_\db)$ be a generalized extension of $A$ of depth $\ell-1$ representing $\epsilon\in S^\ast_{\ell-1}(A)$. Fixing $m,n$ in the eligible range, we need to show that the extension \eqref{eq19} is totally nonsplit. The data of the generalized extension $(M_\db)$ allows us to identify $Gr^WM_{m,n}\cong \bigoplus_{m+1\leq r\leq n}A_r$ (see \cite[Lemma 3.5.1(b)]{Es1}). Let
\[(\mathcal{E}_r)_{1\leq r\leq k-1}\in \bigoplus_r Ext^1_\bT(A_{r+1},A_r)\] 
be the truncation of $(M_\db)$ to depth one (i.e., $\mathcal{E}_r$ is given by \eqref{eq34}). By definition of $S^\ast_{\ell-1}(A)$, this truncation belongs to $\bigm(\bigoplus_r Ext^1_\bT(A_{r+1},A_r)\bigm)^\ast$. That is, considered as an element of $Ext^1_\bT(\mathbbm{1}, V)$, the extension $(\mathcal{E}_r)_{1\leq r\leq k-1}$ is totally disjoint from $\langle A\rangle^{\otimes}$. On recalling the definition of total disjointness (Definition \ref{def: a totally disjoint extension}), we see easily that this forces the extension
\[
(\mathcal{E}_r)_{m< r<n} \in Ext^1_\bT \bigm(\mathbbm{1}, \bigoplus_{m< r< n} \inHom(A_{r+1},A_r)\bigm)
\]
to be totally disjoint from $\langle A\rangle^{\otimes}$. A fortiori, the latter extension is totally disjoint from $\langle Gr^WM_{m,n}\rangle^{\otimes}$. By \cite[Lemma 3.7.3(b)]{Es1} (applied to the generalized extension of $\bigoplus_{m+1\leq r\leq n}A_r$ of full depth obtained by cropping $(M_\db)$ to the part to the left and below $M_{m,n}$, see diagram \eqref{eq18}), the extension $(\mathcal{E}_r)_{m< r<n}$ above coincides with the extension
\[
(W_{p_{r+1}}M_{m,n}/W_{p_{r-1}}M_{m,n})_{m< r< n} \in Ext^1_\bT \bigm(\mathbbm{1}, \bigoplus_{m< r< n} \inHom(A_{r+1},A_r)\bigm),
\]
where we have used the identification $Gr^WM_{m,n}\cong \bigoplus_{m+1\leq r\leq n}A_r$ given by the data of our generalized extension to identify $Gr^W_{p_{r}}M_{m,n}$ with $A_r$. In view of Lemma \ref{lem: totally disjoint invariant under the action of Aut(A)}(a), it follows that the element
\[
(W_{p_{r+1}}M_{m,n}/W_{p_{r-1}}M_{m,n})_{m< r< n} \in Ext^1_\bT\bigm(\mathbbm{1}, \bigoplus_{m< r< n} \inHom(Gr^W_{p_{r+1}}M_{m,n},Gr^W_{p_{r}}M_{m,n})\bigm)
\]
is totally disjoint from $\langle Gr^WM_{m,n}\rangle^{\otimes}$. Thus criterion (iii) of Proposition \ref{prop: maximality criteria 1}(a) holds for $M_{m,n}$, so that $\fu(M_{m,n})=W_{-1}\inEnd(M_{m,n})$. By Lemma 5.2.1 of \cite{Es1} (see also Remark 5.2.4 therein), this implies that every extension coming from the weight filtration of $M_{m,n}$, and in particular the one in \eqref{eq18}, is totally nonsplit. \hfill \qedsymbol{}

\end{document}